\documentclass[a4paper,12pt]{article}
\usepackage{amssymb}
\usepackage{amsthm}
\usepackage{amsmath}
\usepackage[mathscr]{eucal}
\usepackage{fullpage}
\newtheorem{thm}{Theorem}[section]
\newtheorem{lem}[thm]{Lemma}
\newtheorem{cor}[thm]{Corollary}
\newtheorem{prop}[thm]{Proposition}
\theoremstyle{remark}
\newtheorem{rem}[thm]{Remark}
\newtheorem{exmp}[thm]{Example}
\theoremstyle{definition}
\newtheorem{defn}[thm]{Definition}
\title{A characterization of finitely generated reflection subgroups of Coxeter groups orthogonal to a reflection}
\author{Koji Nuida}
\date{Research Center for Information Security (RCIS), National Institute of Advanced Industrial Science and Technology (AIST) \\ %
k.nuida[at]aist.go.jp %
}
\begin{document}
%\markright{Reflection subgroups of Coxeter groups}
%\pagestyle{myheadings}

\maketitle
\begin{abstract}
Given a reflection $r$ in a Coxeter group $W$ (possibly of infinite rank), we consider the subgroup of $W$ generated by the reflections in $W$ having ($-1$)-eigenvectors orthogonal to the ($-1$)-eigenvector of $r$.
In this paper, we determine completely when this subgroup is finitely generated, by using preceding results on centralizers of reflections.
This result provides a new and naturally constructed example of a family of non-finitely generated subgroups of finitely generated groups.
\end{abstract}

\section{Introduction}
\label{sec:intro}

A property for abstract groups is called (e.g., in Section 2.3 of \cite{Rob}) a group-theoretical property if the trivial group possesses the property and any group isomorphic to a group with this property also possesses the property.
Whether or not a given group-theoretical property is preserved by taking subgroups is an important problem with long history in combinatorial group theory.
For affirmative directions, it is a trivial example that the properties of being finite and of being Abelian belong to this class of properties.
Another example which is classical and far from being trivial is that freeness of groups belongs to this class, by virtue of the Nielsen--Schreier Theorem.
On the other hand, of course not every group-theoretical property belongs to this class; a classical but not so obvious example is that a subgroup of a finitely generated group is not necessarily finitely generated.
A famous example of this phenomenon comes from the kernel of the group homomorphism from the free group of rank two, with generators $a$ and $b$, to an infinite cyclic group defined by $a \mapsto a$, $b \mapsto 1$; now the kernel has countably infinite free generators $a^{-i}ba^i$, $i \in \mathbb{Z}$ (see e.g., Example 3(3) of \cite{Bau}).
In this paper, we present a new and naturally constructed example of a family of non-finitely generated subgroups of finitely generated groups, by using the theory of Coxeter groups.
(We notice that another example of such a family based on Coxeter groups was provided by S.~R.~Gal \cite{Gal}; see \cite[Proposition 2.1 and Remark 2.2]{Gal}.)

Let $(W,S)$ be an arbitrary Coxeter system, that is, $W$ is a Coxeter group and $S$ is the distinguished generating set of $W$.
We call an element $w \in W$ a reflection in $W$ if $w$ is conjugate in $W$ to some element of $S$.
In the canonical orthogonal reflection representation of $W$ (over $\mathbb{R}$), an element of $W$ is a reflection in the above sense precisely when it is a reflection in geometrical sense.
Subgroups of Coxeter groups generated by some reflections, referred to as reflection subgroups, are one of the most important kinds of subgroups of Coxeter groups; for example, a theorem of V.~V.~Deodhar \cite{Deo_refsub} or of M.~Dyer \cite[Theorem 3.3]{Dye} states that any reflection subgroup of a Coxeter group is also a Coxeter group (with an appropriately chosen generating set).

In this paper, we deal with a reflection subgroup $W^{\perp r}$ of a Coxeter group $W$ associated to a reflection $r$ in $W$, which is generated by all the reflections in $W$, except $r$ itself, that commute with $r$.
Equivalently, $W^{\perp r}$ is generated by the reflections in $W$ having ($-1$)-eigenvectors which are orthogonal to that of $r$.
The main result of this paper gives a characterization of the cases in which $W^{\perp r}$ is finitely generated (Theorem \ref{thm:condition_for_finite_generation_with_cycle} and Theorem \ref{thm:condition_for_finite_generation_without_cycle}).
Our result shows that the cases in which $W$ is finitely generated but $W^{\perp r}$ is not finitely generated occupies a fairly large part of all the possible cases; hence even such a naturally introduced subgroup of a Coxeter group provides a large family of concrete examples for an important and interesting phenomenon in combinatorial group theory mentioned in the first paragraph.

For the above reflection subgroups $W^{\perp r}$, it holds by definition that $W^{\perp r}$ is a subgroup of the centralizer $Z_W(r)$ of the reflection $r$ in $W$.
The structure of $Z_W(r)$ was studied first by B.~Brink \cite{Bri} in a general setting, and then by R.~E.~Borcherds \cite{Bor}, by Brink and R.~B.~Howlett \cite{Bri-How_normal} and by the author of this paper \cite{Nui_centra} in general and more extended settings.
(We notice that a partial result was presented earlier by Howlett \cite{How}.)
Brink \cite{Bri} gave a semidirect product decomposition of $Z_W(r)$ and determined the structures of some, but not all, factors of $Z_W(r)$.
In particular, Brink's result is not enough to determine the structure of the subgroup $W^{\perp r}$, or even the cardinality of the generating set of $W^{\perp r}$.
On the other hand, the general results in the latter three papers \cite{Bor,Bri-How_normal,Nui_centra} (applied to this special case) determine the entire structure of $Z_W(r)$.
In this paper, we study the detailed structure of $W^{\perp r}$ by using the result of \cite{Nui_centra}.

This paper is organized as follows.
In Section \ref{sec:preliminary}, we summarize some basic definitions and facts about graphs, groupoids (Section \ref{sec:preliminary_graph}) and Coxeter groups (Section \ref{sec:preliminary_Coxetergroup}).
In Section \ref{sec:centralizer}, we describe the preceding results on the centralizers of reflections relevant to our argument in this paper, mainly based on the papers by Brink \cite{Bri} and the author \cite{Nui_centra}.
In Section \ref{sec:characterization}, we first give the statement of the main theorem of this paper as well as its corollaries applied to some special cases (Section \ref{sec:characterization_theorems}), and then give a proof of the main theorem (Section \ref{sec:main_theorem_proof}).

\paragraph*{Acknowledgement.}
The author would like to express his deep gratitude to Professor Itaru Terada and to Professor Kazuhiko Koike for their precious advice and encouragement.
A large part of this work was done during the period when the author was supported by JSPS Research Fellowship (No.\ 16-10825).

\section{Preliminaries}
\label{sec:preliminary}

\subsection{Graphs}
\label{sec:preliminary_graph}

In this paper, any graph $G = (V(G),E(G))$ is simple (i.e., having neither loops nor multiple edges) and unoriented unless otherwise noticed, while it may have infinitely many vertices.
A path in $G$ is denoted by a finite sequence $P = (x_n,\dots,x_1,x_0)$ of its vertices, where the vertices are ordered \emph{from right to left} (i.e., $(x_1,x_0)$ is the first edge in $P$) for a technical reason.
Let $V(P)$ denote the set $\{x_0,x_1,\dots,x_n\}$ of the vertices in the path $P$.
A path $(x_n,\dots,x_1,x_0)$ is said to be \emph{closed} if $x_n = x_0$; and \emph{trivial} if $n = 0$.
A path $(x_0,x_1,\dots,x_n)$ is said to be the \emph{inverse} of $P = (x_n,\dots,x_1,x_0)$ and is denoted by $P^{-1}$.
For two paths $P_1 = (x_n,\dots,x_1,x_0)$ and $P_2 = (y_m,\dots,y_1,y_0)$ with $x_n = y_0$, let $P_2 P_1 := (y_m,\dots,y_1,x_n,\dots,x_1,x_0)$ be the \emph{concatenation} of these two paths.

We say that a path $P = (x_n,\dots,x_1,x_0)$ in $G$ is \emph{simple} if for every $i,j \in \{0,1,\dots,n\}$ with $i \neq j$, we have $x_i \neq x_j$ unless $\{i,j\} = \{0,n\}$.
We say that $P$ is \emph{reduced} if $x_i \neq x_{i+2}$ for every $0 \leq i \leq n-2$; and \emph{cyclically reduced} if $x_0 = x_n$ and $x_i \neq x_j$ for every $i,j \in \{0,1,\dots,n\}$ with $|j-i| \equiv 2 \pmod{n}$.
When $x_n \neq x_0$, we say that $P$ is \emph{full} if for each $0 \leq i \leq j \leq n$ with $j - i \geq 2$, the vertices $x_i$ and $x_j$ are not adjacent in $G$.
On the other hand, when $x_n = x_0$, we say that $P$ is \emph{full} if for each $i,j \in \{0,1,\dots,n\}$ with $|j-i| \not\equiv 0,1 \pmod{n}$, the vertices $x_i$ and $x_j$ are not adjacent in $G$.
We say that a closed path is a \emph{cycle} if it is simple and cyclically reduced.
A \emph{cyclic shift} of a closed path $C = (x_n,\cdots,x_1,x_0)$ is a closed path of the form $C^{\to k} := (x_k,\cdots,x_1,x_0 = x_n,x_{n-1},\dots,x_{k+1},x_k)$ with $0 \leq k \leq n-1$.
We notice the following well-known fact:
\begin{lem}
\label{lem:cycle_and_full_cycle}
Let $(x_n,\dots,x_1,x_0)$ be a cycle in a graph $G$ which is not full.
Then there exist two indices $i$ and $j$ with $0 \leq i \leq j \leq n-1$ and $2 \leq j - i \leq n-2$ satisfying that $(x_i,x_j,\dots,x_{i+1},x_i)$ is a full cycle in $G$ and $(x_n,\dots,x_{j+1},x_j,x_i,\dots,x_1,x_0)$ is a cycle in $G$.
\end{lem}

For a subset $I \subseteq V(G)$, let $G_I$ denote the restriction of $G$ to a vertex set $I$, which is the subgraph of $G$ with vertex set $I$ and edge set consisting of all edges of $G$ whose two vertices are both in $I$.
For a vertex $x \in V(G)$, we define $G_{\sim x}$ to be the connected component of $G$ containing $x$.
\begin{defn}
\label{defn:pre_core}
For a graph $G$, we say that a subset $I \subseteq V(G)$ is a \emph{pre-core} of $G$ if $G_I$ is connected and contains the vertex set of every nontrivial full cycle in $G$.
In this case, we also call $G_I$ a \emph{pre-core} of $G$.
\end{defn}
It follows from Lemma \ref{lem:cycle_and_full_cycle} and induction on the length of a cycle that any pre-core of a graph $G$ contains the vertex set of every (not necessary full) nontrivial cycle in $G$.
Now the following property holds:
\begin{lem}
\label{lem:pre_core_is_convex}
Let $I \subseteq V(G)$ be a pre-core of a graph $G$.
Then every reduced (but not necessarily cyclically reduced) path $P = (x_n,\dots,x_1,x_0)$ in $G$ with $x_0,x_n \in I$ is contained in $G_I$.
\end{lem}
\begin{proof}
We use induction on $n$, with the case $n \leq 1$ being trivial.
We suppose that $n \geq 2$.
First, we consider the case that $x_i = x_j$ for two indices $i < j$ with $(i,j) \neq (0,n)$.
Take such a pair $(i,j)$ with $j - i$ being minimal (hence $j - i \geq 3$, since $P$ is reduced).
Then $(x_j,x_{j-1},\dots,x_{i+1},x_i)$ is a nontrivial cycle in $G$, therefore (by the above remark) it is contained in the pre-core $G_I$.
In particular, $x_i,x_j \in I$.
Now the induction hypothesis implies that both $(x_n,\dots,x_{j+1},x_j)$ and $(x_i,\dots,x_1,x_0)$ are also contained in $G_I$, therefore the claim holds in this case.

Secondly, we consider the other case that $P$ is simple.
By the induction hypothesis, the claim follows if $x_i \in I$ for some $1 \leq i \leq n-1$.
From now, assume contrary that $x_i \not\in I$ for every $1 \leq i \leq n-1$ and deduce a contradiction.
Since the pre-core $I$ is connected, there exists a path $P' = (y_m,\dots,y_1,y_0)$ in $G_I$ with $y_0 = x_0$ and $y_m = x_n$.
By choosing such a path with shortest length, we may assume without loss of generality that $P'$ is simple and reduced when $x_n \neq x_0$; while $P'$ is trivial when $x_n = x_0$.
In any case, since both $P$ and $P'$ are simple and reduced, and $x_i \not\in I$ for every $1 \leq i \leq n-1$ and $y_j \in I$ for every $0 \leq j \leq m$, it follows that $P'{}^{-1} P$ is a nontrivial cycle in $G$.
However, it cannot be contained in $G_I$, since $x_1 \not\in I$.
This is a contradiction.
Hence the proof of Lemma \ref{lem:pre_core_is_convex} is concluded.
\end{proof}

By virtue of Lemma \ref{lem:pre_core_is_convex}, we have the following result:
\begin{lem}
\label{lem:core_exists}
Let $G$ be a graph containing a nontrivial cycle.
Then the intersection of any non-empty family of pre-cores of $G$ is also a pre-core of $G$.
Hence, if all the nontrivial full cycles in such a graph $G$ are contained in the same connected component of $G$, then there exists a unique minimal pre-core of $G$ (which is referred to as the \emph{core} of $G$).
\end{lem}
\begin{proof}
The second part of the claim follows from the first part, since the connected component of $G$ specified in the statement is a pre-core of $G$.
For the first claim, the nontrivial part is that the intersection $I = \bigcap_{\lambda} I_\lambda$ of a non-empty family $\{I_\lambda\}_{\lambda}$ of pre-cores of $G$ is connected.
Choose an index $\lambda = \lambda_0$.
Then for any $y,z \in I$, we have $y,z \in I_{\lambda_0}$, and since $I_{\lambda_0}$ is connected, there exists a reduced path in $G_{I_{\lambda_0}}$ from $y$ to $z$.
Now by Lemma \ref{lem:pre_core_is_convex}, the reduced path is contained in every $G_{I_{\lambda}}$, hence in $G_I$.
This implies that $G_I$ is connected, as desired.
\end{proof}

A \emph{groupoid} is a small category whose morphisms are all invertible; that is, a family of sets $\mathcal{G} = \{\mathcal{G}_{x,y}\}_{x,y \in V(\mathcal{G})}$ with some index set (referred to as \emph{vertex set}) $V(\mathcal{G})$ endowed with (I) multiplications $\mathcal{G}_{x,y} \times \mathcal{G}_{y,z} \to \mathcal{G}_{x,z}$ satisfying the associativity law, (II) an identity element $1=1_x$ in each $\mathcal{G}_{x,x}$, and (III) an inverse $g^{-1} \in \mathcal{G}_{y,x}$ for every $g \in \mathcal{G}_{x,y}$.
Each $\mathcal{G}_{x,x}$ forms a group, which is called a \emph{vertex group} of $\mathcal{G}$.
A \emph{homomorphism} between groupoids is a covariant functor between them regarded as categories.
Notions such as \emph{isomorphisms} and \emph{subgroupoids} are defined as usual.
Now the \emph{fundamental groupoid} $\pi_1(G;\ast,\ast) = \{\pi_1(G;y,x)\}_{x,y}$ of a graph $G$ is defined in the following manner.
The vertex set is chosen to be $V(G)$.
For $x,y \in V(G)$, define $\pi_1(G;y,x)$ to be the quotient set of the set of all paths $(z_n,\dots,z_1,z_0)$ in $G$ from $x = z_0$ to $y = z_n$, with equivalence relation $\sim$ induced by the property that $(z_n,\dots,z_1,z_0) \sim (z_n,\dots,z_{i+3},z_i,\dots,z_1,z_0)$ provided $z_{i+2} = z_i$.
Let $[P]$ denote the equivalence class of a path $P$.
The multiplication $\ast$ in $\pi_1(G;\ast,\ast)$ is given by $[P_1] \ast [P_2] := [P_1P_2]$, which is well-defined.
Then $[P^{-1}]$ is the inverse $[P]^{-1}$ of $[P]$ in $\pi_1(G;\ast,\ast)$ for each path $P$ in $G$.
It is known that each element of $\pi_1(G;\ast,\ast)$ is represented by a unique reduced path in $G$, and $\pi_1(G;\ast,\ast)$ is freely generated by (equivalence classes of) the edges of $G$ (regarded as paths of length $1$).
For each $x \in V(G)$, the vertex group $\pi_1(G;x) := \pi_1(G;x,x)$ is the \emph{fundamental group} of $G$ at $x$, which is known to be a free group.
Note that, if $H$ is a subgraph of $G$, then $\pi_1(H;\ast,\ast)$ is naturally embedded into $\pi_1(G;\ast,\ast)$ as a subgroupoid.
\begin{lem}
\label{lem:fundamental_group_for_subgraph}
Let $G$ be graph.
If $I \subseteq V(G)$ is a pre-core of $G$, then $\pi_1(G_I;x,y) = \pi_1(G;x,y)$ for every $x,y \in I$.
Conversely, if $x \in I \subseteq V(G)$ and $\pi_1(G_I;x) = \pi_1(G;x)$, then $(G_I)_{\sim x}$ contains every nontrivial cycle in $G_{\sim x}$.
\end{lem}
\begin{proof}
The first part of the claim follows from Lemma \ref{lem:pre_core_is_convex}, since each element of $\pi_1(G;x,y)$ is represented by a reduced path between $x$ and $y$.
For the second part, let $C$ be a nontrivial cycle in $G_{\sim x}$.
Take a reduced path $P$ with minimal length in $G_{\sim x}$ from $x$ to some vertex of $C$.
Then by taking an appropriate cyclic shift $C'$ of $C$, it holds that $P^{-1}C'P$ is a reduced closed path in $G$ starting at $x$.
Now we have $[P^{-1}C'P] \in \pi_1(G;x) = \pi_1(G_I;x)$ by the hypothesis.
Since $P^{-1}C'P$ is reduced as above, the uniqueness of representation of each element of the fundamental groupoid of a graph by a reduced path implies that the path $P^{-1}C'P$ is contained in $G_I$, hence in $(G_I)_{\sim x}$.
Therefore $C$ is contained in $(G_I)_{\sim x}$, as desired.
\end{proof}

\subsection{Coxeter groups}
\label{sec:preliminary_Coxetergroup}

A pair $(W,S)$ of a group $W$ and its (possibly infinite) generating set $S$ is called a \emph{Coxeter system} if $W$ admits the following presentation
\begin{displaymath}
W=\langle S \mid (st)^{m(s,t)}=1 \mbox{ for all } s,t \in S \mbox{ with } m(s,t)<\infty \rangle \enspace,
\end{displaymath}
where $m \colon (s,t) \mapsto m(s,t) \in \{1,2,\dots\} \cup \{\infty\}$ is a symmetric mapping in $s,t \in S$ with the property that we have $m(s,t)=1$ if and only if $s=t$.
A group $W$ is called a \emph{Coxeter group} if $(W,S)$ is a Coxeter system for some $S \subseteq W$.
An \emph{isomorphism of Coxeter systems} from $(W,S)$ to $(W',S')$ is a group isomorphism $W \to W'$ that maps $S$ onto $S'$.
It is known that $m(s,t)$ is precisely the order of the element $st$ in $W$, therefore the system $(W,S)$ determines uniquely the mapping $m$ and hence the \emph{Coxeter graph} $\Gamma$, which is a simple undirected graph with vertex set $S$ in which two vertices $s,t \in S$ are joined by an edge with label $m(s,t)$ if and only if $m(s,t) \geq 3$ (by usual convention, the label is omitted when $m(s,t)=3$).
See e.g., the book \cite{Hum} for fundamental properties for Coxeter groups which may be implicit in this paper.

For $w \in W$, let $\ell(w)$ denote the length of $w$ (with respect to $S$), which is the smallest integer $n \geq 0$ with the property that $w$ is expressed by the product of $n$ elements of $S$.
An expression $w = s_1s_2 \cdots s_n$ of $w \in W$ as a product of elements $s_i$ of $S$ with $n = \ell(w)$ is called a \emph{reduced expression} of $w$.
It is known that the subset $\{s_1,s_2,\dots,s_n\} \subseteq S$ is unique and independent of the choice of a reduced expression $w = s_1s_2 \cdots s_n$ of $w$; the set $\{s_1,s_2,\dots,s_n\}$ is called the \emph{support} of $w$ and denoted by $\mathrm{supp}(w)$.
We notice the following properties:
\begin{lem}
\label{lem:support_is_sum}
Let $w_i \in W$ for $1 \leq i \leq k$, $w = w_1w_2 \cdots w_k$, and suppose that $\ell(w) = \sum_{i=1}^{k}\ell(w_i)$.
Then $\mathrm{supp}(w) = \bigcup_{i=1}^{k}\mathrm{supp}(w_i)$.
\end{lem}
\begin{lem}
\label{lem:finitely_generated_iff_finite_rank}
A Coxeter group $W$ is finitely generated if and only if the generating set $S$ is finite.
\end{lem}
\begin{proof}
This follows from the fact that the union of the supports of finitely many elements of $W$ is finite and no proper subset of $S$ can generate $W$.
\end{proof}

For $I \subseteq S$, the subgroup $W_I := \langle I \rangle$ of $W$ generated by $I$ is called a \emph{parabolic subgroup}.
It is well known that $(W_I,I)$ is a Coxeter system with Coxeter graph $\Gamma_I$ and length function $\ell_I$ given by $\ell_I(w) = \ell(w)$ for $w \in W_I$, and we have $\bigcap_{\lambda} W_{J_{\lambda}} = W_J$ (where $J = \bigcap_{\lambda} J_{\lambda}$) for any family $\{J_{\lambda}\}_{\lambda}$ of subsets of $S$.
If $I$ is the vertex set of a connected component of $\Gamma$, then $W_I$ and $I$ are called an \emph{irreducible component} of $W$ (or of $(W,S)$, if we emphasize the generating set $S$) and of $S$, respectively.
Now $W$ is the (restricted) direct product of its irreducible components.
If $\Gamma$ is connected, then $(W,S)$, $W$ and $S$ are called \emph{irreducible}.

We say that a finite subset $I$ of $S$, or the corresponding Coxeter graph $\Gamma_I$, is \emph{of type $\widetilde{A}_{n-1}$} if $n = |I| \geq 3$ and there exists a labelling $I = \{s_1,s_2,\dots,s_n\}$ of elements of $I$ with the property that $m(s_i,s_j) = 3$ for every $i,j \in \{1,2,\dots,n\}$ with $|j-i| \equiv 1 \pmod{n}$ and $m(s_i,s_j) = 2$ for every distinct $i,j \in \{1,2,\dots,n\}$ with $|j-i| \not\equiv 1 \pmod{n}$.

Let $V$ denote the \emph{geometric representation space} (over $\mathbb{R}$), which is an $\mathbb{R}$-linear space equipped with a basis $\Pi = \{\alpha_s \mid s \in S\}$ and a $W$-invariant symmetric bilinear form $\langle \,,\, \rangle$ determined by
\begin{displaymath}
\langle \alpha_s, \alpha_t \rangle =
\begin{cases}
-\cos(\pi / m(s,t)) & \mbox{if } m(s,t) < \infty \enspace; \\
-1 & \mbox{if } m(s,t) = \infty \enspace,
\end{cases}
\end{displaymath}
where $W$ acts faithfully on $V$ by $s \cdot v = v - 2\langle \alpha_s,v \rangle \alpha_s$ for $s \in S$ and $v \in V$.
Then the \emph{root system} $\Phi=W \cdot \Pi$ consists of unit vectors with respect to the bilinear form $\langle \,,\, \rangle$, and $\Phi$ is the disjoint union of $\Phi^+ := \Phi \cap \mathbb{R}_{\geq 0}\Pi$ and $\Phi^- := -\Phi^+$ where $\mathbb{R}_{\geq 0}\Pi$ signifies the set of nonnegative linear combinations of elements of $\Pi$.
Elements of $\Phi$, $\Phi^+$, and $\Phi^-$ are called \emph{roots}, \emph{positive roots}, and \emph{negative roots}, respectively.
We write $\Psi^+ := \Psi \cap \Phi^+$ for any subset $\Psi \subseteq \Phi$.
For an element $v = \sum_{s \in S} c_s \alpha_s \in V$, define the \emph{support} of $v$, denoted by $\mathrm{supp}(v)$, by 
\begin{displaymath}
\mathrm{supp}(v) := \{s \in S \mid c_s \neq 0\} \enspace.
\end{displaymath}
For a subset $I \subseteq S$, let $V_I$ denote the subspace of $V$ spanned by the set $\Pi_I := \{\alpha_s \mid s \in I\} \subseteq \Pi$, and put $\Phi_I := \Phi \cap V_I$.
It is well known that $\Phi_I$ coincides with the root system $W_I \cdot \Pi_I$ of $(W_I,I)$ (see e.g., \cite[Lemma 4]{Fra-How_nearly}).

For a root $\gamma=w \cdot \alpha_s \in \Phi$, let $s_\gamma$ denote the \emph{reflection} $wsw^{-1} \in W$ along the root $\gamma$, acting on $V$ by $s_\gamma \cdot v=v-2\langle \gamma, v \rangle \gamma$ for $v \in V$.
An element $w$ of $W$ is a reflection if and only if $w$ is conjugate in $W$ to an element of $S$.
A subgroup of $W$ generated by some reflections is called a \emph{reflection subgroup}.
The following lemma will be used later:
\begin{lem}
[{\cite[Lemma 2.7]{Nui_ext}}]
\label{lem:wdoesnotfixgamma}
Let $1 \neq w \in W$, $I=\mathrm{supp}(w)$, $\gamma \in \Phi^+$ and $J=\mathrm{supp}(\gamma)$.
Suppose that $I \cap J = \emptyset$ and $J$ is adjacent to $I$ in the Coxeter graph $\Gamma$.
Then we have $w \cdot \gamma \in \Phi_{I \cup J}^+ \smallsetminus \Phi_J$, hence $w \cdot \gamma \neq \gamma$.
\end{lem}

We introduce a special subgraph of a Coxeter graph, which plays an important role in our argument below:
\begin{defn}
\label{defn:odd_Coxeter_graph}
We define the \emph{odd Coxeter graph} $\Gamma^{\mathrm{odd}}$ to be the subgraph of $\Gamma$ obtained by removing all edges labelled by an even number or $\infty$.
\end{defn}

\section{Note on centralizers of reflections}
\label{sec:centralizer}

The structure of the centralizer $Z_W(r)$ of a reflection $r$ in a Coxeter group $W$ has been studied in many preceding papers \cite{Bor,Bri,Bri-How_normal,How,Nui_centra}, in some generalized or restricted settings.
We prepare some notations.
Let $I \subseteq S$, and suppose that $r = s_{\gamma}$ is a reflection in $W_I$ with $\gamma \in \Phi_I$.
We define
\begin{displaymath}
\Phi_I^{\perp r} := \{\beta \in \Phi_I \mid \langle \beta,\gamma \rangle = 0\}
\end{displaymath}
(which is well-defined, since two reflections $s_{\gamma_1}$ and $s_{\gamma_2}$ coincide if and only if $\gamma_2 = \pm \gamma_1$).
Let $W_I^{\perp r}$ denote the subgroup of $W_I$ generated by the reflections $s_{\beta}$ with $\beta \in \Phi_I^{\perp r}$.
We often omit the subscripts \lq $I$' in these notations when $I = S$.
Now if $\gamma = w \cdot \alpha_s$ with $w \in W_I$ and $s \in I$, then the conjugation action by $w$ defines a group isomorphism $W_I^{\perp s} \overset{\sim}{\to} W_I^{\perp r}$, hence $W_I^{\perp r}$ is finitely generated if and only if $W_I^{\perp s}$ is finitely generated.

In what follows, we consider the case $r \in I \subseteq S$.
Let $\Pi^{I,r}$ be the \lq\lq simple system'' in \lq\lq root system'' $\Phi_I^{\perp r}$; that is, $\Pi^{I,r}$ consists of the roots in $(\Phi_I^{\perp r})^+$ that cannot be a positive linear combination of two or more distinct roots in $(\Phi_I^{\perp r})^+$.
We define
\begin{displaymath}
R^{I,r} := \{s_\gamma \mid \gamma \in \Pi^{I,r}\} \enspace.
\end{displaymath}
Then $(W_I^{\perp r},R^{I,r})$ is a Coxeter system by a general theorem of V.~V.~Deodhar \cite{Deo_refsub} or of M.~Dyer \cite[Theorem 3.3]{Dye} on reflection subgroups.
We write $\Pi^{I,r}$ as $\Pi^{r}$ and $R^{I,r}$ as $R^{r}$ when $I = S$.
Now let $Y^I$ be the groupoid, with vertex set $I$ and multiplication induced by that of $W$, defined by
\begin{displaymath}
Y_{y,x}^I := \{w \in W_I \mid \alpha_y = w \cdot \alpha_x \mbox{ and } w \cdot (\Phi_I^{\perp x})^+ \subseteq \Phi^+\} \mbox{ for } x,y \in I
\end{displaymath}
(note that $Y^I$ is indeed a groupoid; see \cite{Bri} or \cite[Section 3.1]{Nui_centra}).
Then B.~Brink \cite{Bri} showed (in slightly different notations) that $Z_{W_I}(r)$ admits the following decomposition:
\begin{displaymath}
Z_{W_I}(r) = (\langle r \rangle \times W_I^{\perp r}) \rtimes Y_{r,r}^I = \langle r \rangle \times (W_I^{\perp r} \rtimes Y_{r,r}^I) \enspace.
\end{displaymath}
Moreover, Brink also proved the following result:
\begin{thm}
\label{thm:Brink_isomorphism}
In this setting, there exists a groupoid isomorphism $\pi^I$ from the fundamental groupoid $\pi_1(\Gamma_I^{\mathrm{odd}};\ast,\ast)$ of the odd Coxeter graph $\Gamma_I^{\mathrm{odd}}$ of $(W_I,I)$ (see Definition \ref{defn:odd_Coxeter_graph} for the definition) to $Y^I$ determined by the following conditions: $\pi^I$ maps each vertex $x \in I$ of $\pi_1(\Gamma_I^{\mathrm{odd}};\ast,\ast)$ to the vertex $x$ of $Y^I$; we have
\begin{displaymath}
\pi^I(y,x) = (xy)^{(m-1)/2} \in Y_{y,x}^I
\end{displaymath}
for any distinct $x,y \in I$ with $m := m(x,y)$ being an odd integer; and we have
\begin{displaymath}
\ell(\pi^I(x_n,\dots,x_1,x_0)) = \sum_{i=1}^{n} \ell(\pi^I(x_i,x_{i-1}))
\end{displaymath}
for any reduced path $(x_n,\dots,x_1,x_0)$ in $\Gamma_I^{\mathrm{odd}}$, where for any path $P = (x_n,\dots,x_1,x_0)$ in $\Gamma_I^{\mathrm{odd}}$, we write $\pi^I([P])$ as $\pi^I(x_n,\dots,x_1,x_0)$ for simplicity.
In particular, $\pi^I$ induces a group isomorphism from $\pi_1(\Gamma_I^{\mathrm{odd}};x)$ to $Y_{x,x}^I$ for each $x \in I$.
\end{thm}
We often omit the superscripts \lq $I$' in the above notations when $I = S$.
Now we present a corollary of Theorem \ref{thm:Brink_isomorphism} which will be used in our argument below:
\begin{cor}
\label{cor:support_of_reduced_path}
Let $I \subseteq S$, and let $P$ be a nontrivial reduced path in $\Gamma_I^{\mathrm{odd}}$.
Then $\mathrm{supp}(\pi^I([P])) = V(P)$.
\end{cor}
\begin{proof}
Since we have $\mathrm{supp}(\pi^I(y,x)) = \{x,y\}$ for any edge $(y,x)$ in $\Gamma_I^{\mathrm{odd}}$ from $x$ to $y$, the claim follows from Theorem \ref{thm:Brink_isomorphism} and Lemma \ref{lem:support_is_sum}.
\end{proof}

The structure of the Coxeter system $(W_I^{\perp r},R^{I,r})$ has been described by a result of the author \cite{Nui_centra} (or a result of Brink and R.~B.~Howlett \cite{Bri-How_normal} in different notations).
For each $r \in I$, let $\mathrm{Odd}_I(r)$ denote the vertex set of $(\Gamma_I^{\mathrm{odd}})_{\sim r}$ (see Section \ref{sec:preliminary_graph} for the notation).
Lemma 4.1 of \cite{Nui_centra} implies the following property:
\begin{lem}
[{see \cite[Lemma 4.1]{Nui_centra}}]
\label{lem:loop_generator}
Let $x$ and $s$ be distinct elements of $S$.
Then we have $\Phi_{\{x,s\}}^{\perp x} \neq \emptyset$ if and only if $m(x,s) < \infty$ and the longest element of the finite Coxeter group $W_{\{x,s\}}$ commutes with $x$.
Moreover, if these two conditions are satisfied, then $(\Phi_{\{x,s\}}^{\perp x})^+$ consists of a unique element, denoted here by $\gamma(x,s)$.
\end{lem}
A direct calculation shows that $x$ and $s$ satisfy the conditions in Lemma \ref{lem:loop_generator} if and only if $m(x,s)$ is an even integer.
Now we define
\begin{displaymath}
\mathcal{E}^I(r) := \{(x,s) \mid x \in \mathrm{Odd}_I(r), s \in I \smallsetminus \{x\} \mbox{ and } m(x,s) \mbox{ is even}\} \enspace,
\end{displaymath}
therefore for each $\xi = (x,s) \in \mathcal{E}^I(r)$, the positive root $\gamma(\xi) = \gamma(x,s)$ as in Lemma \ref{lem:loop_generator} exists.
For each $\xi \in \mathcal{E}^I(r)$, we write $\xi = (\xi^{\circ},\xi_{\dagger})$ with $\xi^{\circ} \in \mathrm{Odd}_I(r)$ and $\xi_{\dagger} \in I \smallsetminus \{\xi^{\circ}\}$, and we call the elements $\xi^{\circ}$ and $\xi_{\dagger}$ the \emph{ship} and \emph{anchor} of $\xi$, respectively.
Now the results of \cite{Nui_centra}, especially Theorem 4.13 combined with Example 4.12, imply the following property (recall the isomorphism $\pi^I$ given in Theorem \ref{thm:Brink_isomorphism}):
\begin{thm}
[{see \cite{Nui_centra}}]
\label{thm:generator_of_orthogonal_reflections}
In this setting, $\Pi^{I,r}$ consists of the roots $\gamma_r^I(c;\xi) = \gamma_r^I(c;\xi^{\circ},\xi_{\dagger})$ with $\xi \in \mathcal{E}^I(r)$ and $c \in \pi_1(\Gamma_I^{\mathrm{odd}};r,\xi^{\circ})$, defined by
\begin{displaymath}
\gamma_r^I(c;\xi) := \pi^I(c) \cdot \gamma(\xi) \enspace.
\end{displaymath}
Hence $R^{I,r}$ consists of the reflections $s_r^I(c;\xi)$ along the roots $\gamma_r^I(c;\xi)$ with $\xi$ and $c$ as above.
\end{thm}
We define
\begin{displaymath}
\mathcal{R}^I(r) := \{(c;\xi) \mid \xi \in \mathcal{E}^I(r) \mbox{ and } c \in \pi_1(\Gamma_I^{\mathrm{odd}};r,\xi^{\circ})\} \enspace.
\end{displaymath}
Then by Theorem \ref{thm:generator_of_orthogonal_reflections}, an element of $\Pi^{I,r}$ (hence of $R^{I,r}$) is determined for each element of $\mathcal{R}^I(r)$.

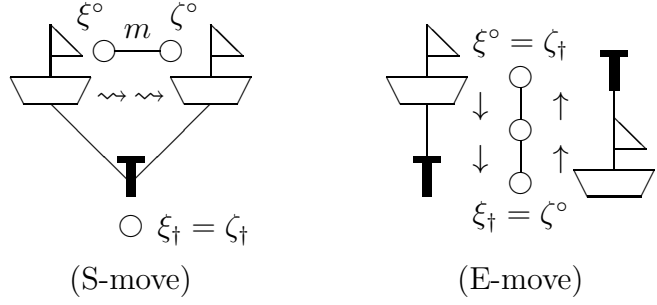
\begin{figure}[hbt]
\centering
\begin{picture}(90,110)(0,-120)
\put(0,-40){\line(1,0){30}}
\put(0,-40){\line(1,-2){5}}
\put(30,-40){\line(-1,-2){5}}
\put(5,-50){\line(1,0){20}}
\put(15,-40){\line(0,1){20}}
\put(15,-20){\line(1,-1){12}}
\put(15,-32){\line(1,0){12}}
\put(35,-30){\circle{8}}
\put(30,-20){\hbox to0pt{\hss$\xi^{\circ}$\hss}}
\put(39,-30){\line(1,0){17}}
\put(60,-30){\circle{8}}
\put(65,-20){\hbox to0pt{\hss$\zeta^{\circ}$\hss}}
\put(47.5,-25){\hbox to0pt{\hss$m$\hss}}
\put(60,-40){\line(1,0){30}}
\put(60,-40){\line(1,-2){5}}
\put(90,-40){\line(-1,-2){5}}
\put(65,-50){\line(1,0){20}}
\put(75,-40){\line(0,1){20}}
\put(75,-20){\line(1,-1){12}}
\put(75,-32){\line(1,0){12}}
\put(15,-50){\line(1,-1){30}}
\put(75,-50){\line(-1,-1){30}}
\multiput(43.5,-70)(.5,0){7}{\line(0,-1){15}}
\multiput(40,-70)(0,-.5){5}{\line(1,0){10}}
\put(45,-96){\circle{8}}
\put(55,-100){$\xi_{\dagger} = \zeta_{\dagger}$}
\put(45,-50){\hbox to0pt{\hss$\leadsto\,\leadsto$\hss}}
\put(45,-120){\hbox to0pt{\hss(S-move)\hss}}
\end{picture}
\qquad\qquad
\begin{picture}(100,110)(0,-120)
\put(0,-40){\line(1,0){30}}
\put(0,-40){\line(1,-2){5}}
\put(30,-40){\line(-1,-2){5}}
\put(5,-50){\line(1,0){20}}
\put(15,-40){\line(0,1){20}}
\put(15,-20){\line(1,-1){12}}
\put(15,-32){\line(1,0){12}}
\put(15,-50){\line(0,-1){30}}
\multiput(13.5,-70)(.5,0){7}{\line(0,-1){15}}
\multiput(10,-70)(0,-.5){5}{\line(1,0){10}}
\put(50,-40){\circle{8}}
\put(50,-30){\hbox to0pt{\hss$\xi^{\circ} = \zeta_{\dagger}$\hss}}
\multiput(50,-44)(0,-20){2}{\line(0,-1){12}}
\multiput(50,-60)(0,-20){2}{\circle{8}}
\put(50,-95){\hbox to0pt{\hss$\xi_{\dagger} = \zeta^{\circ}$\hss}}
\multiput(35,-54)(0,-20){2}{\hbox to0pt{\hss$\downarrow$\hss}}
\multiput(65,-54)(0,-20){2}{\hbox to0pt{\hss$\uparrow$\hss}}
\put(70,-75){\line(1,0){30}}
\put(70,-75){\line(1,-2){5}}
\put(100,-75){\line(-1,-2){5}}
\put(75,-85){\line(1,0){20}}
\put(85,-75){\line(0,1){20}}
\put(85,-55){\line(1,-1){12}}
\put(85,-67){\line(1,0){12}}
\put(85,-55){\line(0,1){10}}
\multiput(83.5,-30)(.5,0){7}{\line(0,-1){15}}
\multiput(80,-30)(0,-.5){5}{\line(1,0){10}}
\put(50,-120){\hbox to0pt{\hss(E-move)\hss}}
\end{picture}
\caption{Picture for S-moves and E-moves}
\label{fig:ship_and_anchor}
\end{figure}

Here we give a remark on an intuition behind the names \lq\lq ship'' and \lq\lq anchor'' for the two components of elements of $\mathcal{E}^I(r)$.
Let $(c;\xi) \in \mathcal{R}^I(r)$.
Suppose that $\xi \neq \zeta \in \mathcal{E}^I(r)$, $\zeta_{\dagger} = \xi_{\dagger}$, $m(\xi^{\circ},\xi_{\dagger}) = 2 = m(\zeta^{\circ},\zeta_{\dagger})$ and $m = m(\xi^{\circ},\zeta^{\circ})$ is an odd integer (hence $\xi^{\circ}$ and $\zeta^{\circ}$ are adjacent in $\Gamma_I^{\mathrm{odd}}$).
In this case, the transformation from $\xi$ to $\zeta$ can be regarded as that the ship of $\xi$ moves from the vertex $\xi^{\circ}$ to the vertex $\zeta^{\circ}$ through the corresponding edge in $\Gamma_I^{\mathrm{odd}}$, while the anchor of $\xi$ at the vertex $\xi_{\dagger}$ is fixed (see the left part of Figure \ref{fig:ship_and_anchor}).
We say that this kind of transformation $\xi \mapsto \zeta$ is an \emph{S-move} in $I$ from $\xi$ to $\zeta$ (\lq S' stands for \lq\lq slide'').
Now we have $\gamma_r^I(c;\xi) = \gamma_r^I(c \ast [q^{-1}];\zeta)$ where $q$ denotes the edge from $\xi^{\circ}$ to $\zeta^{\circ}$.
We also call the transformation $(c;\xi) \mapsto (c \ast [q^{-1}];\zeta)$ an \emph{S-move} in $I$, and we say that the edge $q$, or the element $[q]$ of $\pi_1(\Gamma_I^{\mathrm{odd}};\ast,\ast)$, is \emph{realized} by the S-move.

For the fundamental relations of the Coxeter system $(W_I^{\perp r},R^{I,r})$, roughly speaking, the result of \cite{Nui_centra} implies that all the relations in $(W_I^{\perp r},R^{I,r})$ are induced by those defined within some parabolic subgroups of rank $3$.
Let $J$ be a subset of $I$ with the properties that $|J| = 3$, $J \cap \mathrm{Odd}_I(r) \neq \emptyset$ and $|W_J| < \infty$.
Now for any $r' \in J \cap \mathrm{Odd}_I(r)$, a direct calculation shows that there exists an element $\xi \in \mathcal{E}^I(r)$ with $\xi^{\circ} \in \mathrm{Odd}_J(r')$ and $\xi_{\dagger} \in J$.
Moreover, there exist a unique element $\zeta \in \mathcal{E}^I(r) \setminus \{\xi\}$ with $\zeta^{\circ} \in \mathrm{Odd}_J(r') = \mathrm{Odd}_J(\xi^{\circ})$ and $\zeta_{\dagger} \in J$, and a unique (possibly trivial) reduced path $q$ in $\Gamma_J^{\mathrm{odd}}$ from $\xi^{\circ}$ to $\zeta^{\circ}$.
All the possibilities are listed (up to symmetry) in Table \ref{tab:fundamental_relation}, where we put $J = \{x,y,z\}$ and write $I_2(3) = A_2$ and $I_2(4) = B_2$ for simplicity.
Now the argument in the final part of Section 4.4 of \cite{Nui_centra} implies that we have $\langle \gamma(\zeta),\pi^I([q]) \cdot \gamma(\xi) \rangle = -\cos(\pi/k)$ for an integer $k \geq 1$; the column \lq\lq order'' in Table \ref{tab:fundamental_relation} gives the integer $k$ for each case (cf., Tables 1 and 2 in \cite{Nui_centra}).
In this setting, for each $c \in \pi_1(\Gamma_I^{\mathrm{odd}};r,\xi^{\circ})$, we introduce the symmetric relations $(c;\xi) \overset{k}{\sim}_I (c \ast [q^{-1}];\zeta)$ and $(c \ast [q^{-1}];\zeta) \overset{k}{\sim}_I (c;\xi)$ for two elements of $\mathcal{R}^I(r)$.

\begin{table}[hbt]
\centering
\caption{List for the local relations, where $J = \{x,y,z\}$, $I_2(3) = A_2$, $I_2(4) = B_2$}
\label{tab:fundamental_relation}
\begin{tabular}{|c|c|c|c|c|c|c} \cline{1-6}
type of $W_J$ & $m(x,y)$, $m(y,z)$ & $\xi$ & path $q$ & $\zeta$ & order & \\ \cline{1-6}
$A_1 \times A_1 \times A_1$ & $2$, $2$ & $(x,y)$ & trivial & $(x,z)$ & $2$ & \\ \cline{1-6}
$A_1 \times I_2(m)$ & $2$, $m$ & $(x,y)$ & trivial & $(x,z)$ & $m$ & \\ \cline{3-6}
($m \geq 3$ odd) & & $(y,x)$ & $(z,y)$ & $(z,x)$ & $1$ & (S-move) \\ \cline{1-6}
$A_1 \times I_2(m)$ & $2$, $m$ & $(x,y)$ & trivial & $(x,z)$ & $m$ & \\ \cline{3-6}
($m \geq 4$ even) & & $(y,x)$ & trivial & $(y,z)$ & $2$ & \\ \cline{1-6}
$A_3$ & $3$, $3$ & $(x,z)$ & $(z,y,x)$ & $(z,x)$ & $1$ & (E-move) \\ \cline{1-6}
$B_3$ & $3$, $4$ & $(x,z)$ & $(y,x)$ & $(y,z)$ & $2$ & \\ \cline{3-6}
& & $(z,x)$ & trivial & $(z,y)$ & $4$ \\ \cline{1-6}
$H_3$ & $3$, $5$ & $(x,z)$ & $(z,y,x)$ & $(z,x)$ & $2$ & \\ \cline{1-6}
\end{tabular}
\end{table}

We focus on two of the above-mentioned cases with $k = 1$, the third and sixth cases in Table \ref{tab:fundamental_relation}.
The former case corresponds to an S-move introduced above.
For the latter case, namely $\xi = (x,z)$, $\zeta = (z,x)$, $q = (z,y,x)$ and $\{x,y,z\}$ is of type $A_3$, the transformation from $\xi$ to $\zeta$ can be regarded as that the ship of $\xi$ moves from the vertex $\xi^{\circ}$ to the vertex $\zeta^{\circ}$ through the path $q$ in $\Gamma_I^{\mathrm{odd}}$, and the positions of the ship and the anchor are exchanged during the transformation (see the right part of Figure \ref{fig:ship_and_anchor}).
We say that this kind of transformation $\xi \mapsto \zeta$ is an \emph{E-move} in $I$ from $\xi$ to $\zeta$ (\lq E' stands for \lq\lq exchange'').
Now we have $\gamma_r^I(c;\xi) = \gamma_r^I(c \ast [q^{-1}];\zeta)$ for each $c \in \pi_1(\Gamma_I^{\mathrm{odd}};r,\xi^{\circ})$.
We also call the transformation $(c;\xi) \mapsto (c \ast [q^{-1}];\zeta)$ an \emph{E-move} in $I$, and we say that the edge $q$, or the element $[q]$ of $\pi_1(\Gamma_I^{\mathrm{odd}};\ast,\ast)$, is \emph{realized} by the E-move.

Let the term \emph{move} signify both an S-move and an E-move.
Note that each move is invertible by the definition.
For a sequence of $n$ consecutive moves, where the $i$-th move ($1 \leq i \leq n$) realizes a path $q_i$ (or an element $[q_i]$ of the fundamental groupoid), we say that the concatenation $q := q_n \cdots q_1$ of the paths $q_i$ (or the multiplication $[q]$ of the elements $[q_i]$) is \emph{realized} by the sequence of moves.
Now we define an equivalence relation $\sim_I$ on $\mathcal{R}^I(r)$ in such a way that we have $(c;\xi) \sim_I (d;\zeta)$ if and only if there exists a sequence of moves from $\xi$ to $\zeta$ that realizes an element $[q]$ with $d = c \ast [q]^{-1}$.

Now Theorem 4.14 of \cite{Nui_centra} implies the following property:
\begin{thm}
\label{thm:presentation_of_orthogonal_reflection_subgroup}
Let $(c;\xi),(d;\zeta) \in \mathcal{R}^I(r)$.
Then we have $s_r^I(c;\xi) = s_r^I(d;\zeta)$ if and only if $(c;\xi) \sim_I (d;\zeta)$.
Moreover, for $2 \leq m < \infty$, the product of $s_r^I(c;\xi)$ and $s_r^I(d;\zeta)$ has order $m$ if and only if there exist two elements $(c';\xi')$ and $(d';\zeta')$ of $\mathcal{R}^I(r)$ with the property that $(c;\xi) \sim_I (c';\xi')$, $(d;\zeta) \sim_I (d';\zeta')$ and $(c';\xi') \overset{m}{\sim}_I (d';\zeta')$.
\end{thm}
We often omit the subscripts \lq $I$' and superscripts \lq $I$' in the above notations when $I = S$.

\begin{exmp}
\label{ex:cycle_odd_length}
Let $P = (x_n,\dots,x_1,x_0)$ be a full cycle in $\Gamma^{\mathrm{odd}}$.
Write $\xi_{i,j} = (x_i,x_j)$ and $m_{i,j} = m(x_i,x_j)$ for simplicity, where the indices $i,j$ are taken modulo $n$.
In this example, we consider the case that $n = 2N + 1$ with an integer $N \geq 2$ and $m_{0,k} < \infty$ for some index $2 \leq k \leq n-2$.
Now we have $\xi_{k,0} \in \mathcal{E}^{V(P)}(x_0)$ for this $k$.

We consider the possibilities of moves in $V(P)$ starting at $\xi_{i,j}$, where $\xi_{i,j}$ is an element of $\mathcal{E}^{V(P)}(x_0)$ (hence $i-j \not\equiv 0,\pm 1 \pmod{n}$).
By the shape of $\Gamma_{V(P)}^{\mathrm{odd}}$, for the case $i-j \not\equiv \pm 2 \pmod{n}$, there exist at most two moves in $V(P)$ starting at $\xi_{i,j}$, the S-move $\xi_{i,j} \mapsto \xi_{i-1,j}$ and the S-move $\xi_{i,j} \mapsto \xi_{i+1,j}$.
Moreover, the former S-move $\xi_{i,j} \mapsto \xi_{i-1,j}$ indeed exists if and only if $m_{i,j} = m_{i-1,j} = 2$, while the latter S-move $\xi_{i,j} \mapsto \xi_{i+1,j}$ indeed exists if and only if $m_{i,j} = m_{i+1,j} = 2$.
On the other hand, for the case $i-j \equiv 2 \pmod{n}$, there exist at most two moves in $V(P)$ starting at $\xi_{i,j} = \xi_{j+2,j}$, the E-move $\xi_{j+2,j} \mapsto \xi_{j,j+2}$ in $\{x_j,x_{j+1},x_{j+2}\}$ and the S-move $\xi_{j+2,j} \mapsto \xi_{j+3,j}$.
Moreover, the former E-move $\xi_{j+2,j} \mapsto \xi_{j,j+2}$ indeed exists if and only if $m_{j,j+2} = 2$ and $m_{j,j+1} = m_{j+1,j+2} = 3$, while the latter S-move $\xi_{j+2,j} \mapsto \xi_{j+3,j}$ indeed exists if and only if $m_{j,j+2} = m_{j,j+3} = 2$.
A similar property holds for the remaining case $i-j \equiv -2 \pmod{n}$.

By the above argument, if a non-backtracking sequence of moves in $V(P)$ from $\xi_{k,0}$ to $\xi_{k,0}$ exists, then it is a multiple of the concatenation of the following sequences (recall that $n = 2N+1$) and the conditions specified there are satisfied, where $\xi_{i,j} \to \xi_{i',j}$ means the S-move from $\xi_{i,j}$ to $\xi_{i',j}$ and $\xi_{i,j} \overset{h}{\leftrightarrow} \xi_{j,i}$ means the E-move in $\{x_i,x_j,x_h\}$ from $\xi_{i,j}$ to $\xi_{j,i}$:
\begin{itemize}
\item $\xi_{k,0} \to \xi_{k+1,0} \to \cdots \to \xi_{n-2,0} \overset{n-1}{\leftrightarrow} \xi_{0,n-2}$ (Condition: $m_{0,i} = 2$ for every $k \leq i \leq n-2$, and $m_{n-2,n-1} = m_{n-1,0} = 3$);
\item $\xi_{2j+3,2j+1} \to \xi_{2j+4,2j+1} \to \cdots \to \xi_{2j-1,2j+1} \overset{2j}{\leftrightarrow} \xi_{2j+1,2j-1}$ for $j = N-1,N-2,\dots,1$ (Condition: $m_{2j+1,i} = 2$ for every $i$ with $i - (2j+1) \not\equiv 0,\pm 1 \pmod{n}$, and $m_{2j-1,2j} = m_{2j,2j+1} = 3$);
\item $\xi_{3,1} \to \xi_{4,1} \to \cdots \to \xi_{n-1,1} \overset{0}{\leftrightarrow} \xi_{1,n-1}$ (Condition: $m_{1,i} = 2$ for every $3 \leq i \leq n-1$, and $m_{n-1,0} = m_{0,1} = 3$);
\item $\xi_{2j+2,2j} \to \xi_{2j+3,2j} \to \cdots \to \xi_{2j-2,2j} \overset{2j-1}{\leftrightarrow} \xi_{2j,2j-2}$ for $j = N,N-1,\dots,1$ (Condition: $m_{2j,i} = 2$ for every $i$ with $i - 2j \not\equiv 0,\pm 1 \pmod{n}$, and $m_{2j-2,2j-1} = m_{2j-1,2j} = 3$);
\item $\xi_{2,0} \to \xi_{3,0} \to \cdots \to \xi_{k,0}$ (Condition: $m_{0,i} = 2$ for every $2 \leq i \leq k$).
\end{itemize}
Now the fundamental group $\pi_1(\Gamma_{V(P)}^{\mathrm{odd}};x_k)$ is the infinite cyclic group generated by the equivalence class $[P^{\to k}]$ of a cyclic shift of $P$ (see Section \ref{sec:preliminary_graph} for the notation), while the concatenation of the above sequences of moves realizes the reduced closed path $(P^{\to k})^{n-2}$.
Moreover, all the conditions specified above are satisfied if and only if we have $m_{i,j} = 2$ for every $i,j$ with $i-j \not\equiv 0,\pm 1 \pmod{n}$, and $m_{i,j} = 3$ for every $i,j$ with $i-j \equiv \pm 1 \pmod{n}$; namely, $V(P)$ is of type $\widetilde{A}_{n-1}$.

Summarizing, if $V(P)$ is not of type $\widetilde{A}_{n-1}$, then no nontrivial element of $\pi_1(\Gamma_{V(P)}^{\mathrm{odd}};x_k)$ is realized by a sequence of moves in $V(P)$.
On the other hand, if $V(P)$ is of type $\widetilde{A}_{n-1}$, then an element $[P^{\to k}]^{\ell}$ of $\pi_1(\Gamma_{V(P)}^{\mathrm{odd}};x_k)$ is realized by a sequence of moves in $V(P)$ if and only if $\ell$ is a multiple of $n-2$; therefore we have $s_{x_0}([P^{\to k}]^{\ell_1};\xi_{k,0}) = s_{x_0}([P^{\to k}]^{\ell_2};\xi_{k,0})$ if and only if $\ell_1 \equiv \ell_2 \pmod{n-2}$.
Moreover, in any case, (since now $n \geq 5$) the generator $[P^{\to k}]$ of $\pi_1(\Gamma_{V(P)}^{\mathrm{odd}};x_k)$ cannot be realized by a sequence of moves in $V(P)$.
\end{exmp}
\begin{exmp}
\label{ex:cycle_even_length}
We use the same notations as Example \ref{ex:cycle_odd_length}, but now we consider the case that $n = 2N$ with an integer $N \geq 2$.
We also suppose that $m_{0,k} < \infty$ for some index $2 \leq k \leq n-2$, hence $\xi_{k,0} \in \mathcal{E}^{V(P)}(x_0)$ for this $k$.

Now, by the same argument as Example \ref{ex:cycle_odd_length}, if a non-backtracking sequence of moves in $V(P)$ from $\xi_{k,0}$ to $\xi_{k,0}$ exists, then it is a multiple of the concatenation of the following sequences (recall that $n = 2N$) and the conditions specified there are satisfied:
\begin{itemize}
\item $\xi_{k,0} \to \xi_{k+1,0} \to \cdots \to \xi_{n-2,0} \overset{n-1}{\leftrightarrow} \xi_{0,n-2}$ (Condition: $m_{0,i} = 2$ for every $k \leq i \leq n-2$, and $m_{n-2,n-1} = m_{n-1,0} = 3$);
\item $\xi_{2j+2,2j} \to \xi_{2j+3,2j} \to \cdots \to \xi_{2j-2,2j} \overset{2j-1}{\leftrightarrow} \xi_{2j,2j-2}$ for $j = N-1,N-2,\dots,1$ (Condition: $m_{2j,i} = 2$ for every $i$ with $i - 2j \not\equiv 0,\pm 1 \pmod{n}$, and $m_{2j-2,2j-1} = m_{2j-1,2j} = 3$);
\item $\xi_{2,0} \to \xi_{3,0} \to \cdots \to \xi_{k,0}$ (Condition: $m_{0,i} = 2$ for every $2 \leq i \leq k$).
\end{itemize}
The fundamental group $\pi_1(\Gamma_{V(P)}^{\mathrm{odd}};x_k)$ is the infinite cyclic group generated by $[P^{\to k}]$, while the concatenation of the above sequences of moves realizes the reduced closed path $(P^{\to k})^{N-1}$.

The possibilities of $\Gamma_{V(P)}$ are slightly more complicated than Example \ref{ex:cycle_odd_length}.
Namely, all the conditions specified above are satisfied if and only if we have $m_{i,j} = 2$ for every $i,j$ satisfying that $i-j \not\equiv 0,\pm 1 \pmod{n}$ \emph{and at least one of $i$ and $j$ is even}, and $m_{i,j} = 3$ for every $i,j$ with $i-j \equiv \pm 1 \pmod{n}$.

Summarizing, if $\Gamma_{V(P)}$ does not satisfy the above condition, then no nontrivial element of $\pi_1(\Gamma_{V(P)}^{\mathrm{odd}};x_k)$ is realized by a sequence of moves in $V(P)$.
On the other hand, if $\Gamma_{V(P)}$ satisfies the above condition, then an element $[P^{\to k}]^{\ell}$ of $\pi_1(\Gamma_{V(P)}^{\mathrm{odd}};x_k)$ is realized by a sequence of moves in $V(P)$ if and only if $\ell$ is a multiple of $N-1 = n/2 - 1$; therefore we have $s_{x_0}([P^{\to k}]^{\ell_1};\xi_{k,0}) = s_{x_0}([P^{\to k}]^{\ell_2};\xi_{k,0})$ if and only if $\ell_1 \equiv \ell_2 \pmod{n/2-1}$.
Moreover, in any case (since now $n \geq 4$) the generator $[P^{\to k}]$ of $\pi_1(\Gamma_{V(P)}^{\mathrm{odd}};x_k)$ is realized by a sequence of moves in $V(P)$ if and only if $n = 4$ (hence $k = 2$) and we have $m_{0,1} = m_{1,2} = m_{2,3} = m_{3,0} = 3$, $m_{0,2} = 2$, and $m_{1,3}$ is either an even integer or $\infty$.
\end{exmp}

We give some remarks on relations between the Coxeter systems $(W_I^{\perp r},R^{I,r})$ for different $I$ and $r$:
\begin{rem}
\label{rem:restriction_to_smaller_domain}
Suppose that $r \in I \subseteq J \subseteq S$.
Then we have $\mathrm{Odd}_I(r) \subseteq \mathrm{Odd}_J(r)$, therefore $\mathcal{E}^I(r) \subseteq \mathcal{E}^J(r)$ by the definition.
On the other hand, since $\Gamma_I^{\mathrm{odd}}$ is a subgraph of $\Gamma_J^{\mathrm{odd}}$, the groupoid $\pi_1(\Gamma_I^{\mathrm{odd}};\ast,\ast)$ is naturally embedded into $\pi_1(\Gamma_J^{\mathrm{odd}};\ast,\ast)$, and the isomorphism $\pi^J$ agrees with $\pi^I$ on $\pi_1(\Gamma_I^{\mathrm{odd}};\ast,\ast)$.
This implies that $\mathcal{R}^I(r) \subseteq \mathcal{R}^J(r)$ and we have $\gamma_r^I(c;\xi) = \gamma_r^J(c;\xi)$ and $s_r^I(c;\xi) = s_r^J(c;\xi)$ for each $(c;\xi) \in \mathcal{R}^I(r)$; hence we have $\Pi^{I,r} \subseteq \Pi^{J,r}$ and $R^{I,r} \subseteq R^{J,r}$.
\end{rem}
\begin{rem}
\label{rem:isomorphism_by_conjugation}
Suppose that $r \in I \subseteq S$ and $r' \in \mathrm{Odd}_I(r)$.
Then $\mathcal{E}^I(r) = \mathcal{E}^I(r')$ by definition.
On the other hand, there exists a path $p$ in $\Gamma_I^{\mathrm{odd}}$ from $r$ to $r'$.
Then by Theorem \ref{thm:generator_of_orthogonal_reflections}, the conjugation action by $\pi^I([p]) \in Y_{r',r}^I$ maps an element $s_r^I(c;\xi)$ of $R^{I,r}$ to an element $s_{r'}^I([p] \ast c;\xi)$ of $R^{I,r'}$, which gives rise to an isomorphism of Coxeter systems from $(W_I^{\perp r},R^{I,r})$ to $(W_I^{\perp r'},R^{I,r'})$.
\end{rem}

\section{Characterization of finitely generated $W^{\perp r}$}
\label{sec:characterization}

In this section, we determine a necessary and sufficient condition for the reflection subgroup $W^{\perp r}$ with $r \in S$ introduced in Section \ref{sec:centralizer} to be finitely generated.
This results shows that a finitely generated Coxeter group may frequently have a reflection subgroup which is not finitely generated.

\subsection{Statements of main theorems and special cases}
\label{sec:characterization_theorems}

In this subsection, we present the statements of the main theorems of this paper.
Proofs will be given in the next subsection.
We also give some corollaries of the main theorems applied to some special cases.

We introduce some notations.
For $r \in S$, we define
\begin{displaymath}
\begin{split}
\partial\mathrm{EOdd}(r)
:={}& \{x \in S \smallsetminus \mathrm{Odd}(r) \mid m(x,y) < \infty \mbox{ for some } y \in \mathrm{Odd}(r)\} \\
={}& \{x \in S \smallsetminus \mathrm{Odd}(r) \mid m(x,y) \mbox{ is an even integer for some } y \in \mathrm{Odd}(r)\}
\end{split}
\end{displaymath}
where we write $\mathrm{Odd}(r) = \mathrm{Odd}_S(r)$.
We put
\begin{displaymath}
\mathrm{EOdd}(r) := \mathrm{Odd}(r) \cup \partial\mathrm{EOdd}(r) \enspace.
\end{displaymath}
On the other hand, we define
\begin{displaymath}
\mathrm{Odd}(r)^{\perp x} := \{y \in \mathrm{Odd}(r) \mid m(x,y) = 2\} \mbox{ for } x \in S \enspace.
\end{displaymath}
Note that $\mathrm{Odd}(r)^{\perp x} = \emptyset$ if $x \not\in \mathrm{EOdd}(r)$.

Now we give the statement of the main theorem of this paper, which is divided into the following two cases:
\begin{thm}
\label{thm:condition_for_finite_generation_with_cycle}
Let $r \in S$.
Suppose that $(\Gamma^{\mathrm{odd}})_{\sim r}$ contains a nontrivial cycle.
Let $K \subseteq S$ be the core of $(\Gamma^{\mathrm{odd}})_{\sim r}$ (see Lemma \ref{lem:core_exists} for the terminology).
Then $W^{\perp r}$ is finitely generated if and only if all the following four conditions are satisfied (see above for the notations):
\begin{enumerate}
\item \label{item:thm_condition_for_finite_generation_with_cycle__core}
One of the following three conditions is satisfied:
\begin{enumerate}
\item \label{item:thm_condition_for_finite_generation_with_cycle__core_not_even}
$m(x,y)$ is not an even integer for any $x,y \in \mathrm{Odd}(r)$;
\item \label{item:thm_condition_for_finite_generation_with_cycle__core_A_tilde}
$K$ is of type $\widetilde{A}_{n-1}$ with $n \geq 4$, and we have $x,y \in K$ whenever $x,y \in \mathrm{Odd}(r)$ and $m(x,y)$ is an even integer;
\item \label{item:thm_condition_for_finite_generation_with_cycle__core_bipyramid}
$K$ has at least $4$ elements, and there exist two elements $x_1$ and $x_2$ in $K$ satisfying the following conditions: $m(x_1,x_2) = 2$; $m(x_i,y) = 3$ for each $y \in K \smallsetminus \{x_1,x_2\}$ and $i \in \{1,2\}$; $m(y,z) = \infty$ for every distinct $y,z \in K \smallsetminus \{x_1,x_2\}$; and we have $\{y,z\} = \{x_1,x_2\}$ whenever $y,z \in \mathrm{Odd}(r)$ and $m(y,z)$ is an even integer.
\end{enumerate}
\item \label{item:thm_condition_for_finite_generation_with_cycle__E_finite}
The set $\partial\mathrm{EOdd}(r)$ is finite.
\item \label{item:thm_condition_for_finite_generation_with_cycle__2_or_infinite}
We have $m(x,y) \in \{2,\infty\}$ for every $x \in \mathrm{Odd}(r)$ and $y \in \partial\mathrm{EOdd}(r)$.
\item \label{item:thm_condition_for_finite_generation_with_cycle__non-adjacent_vertices}
For each $x \in \partial\mathrm{EOdd}(r)$, the graph $\Gamma_{\mathrm{Odd}(r)^{\perp x}}^{\mathrm{odd}}$ is connected and contains $K$ (or equivalently, $\mathrm{Odd}(r)^{\perp x}$ is a pre-core of $(\Gamma^{\mathrm{odd}})_{\sim r}$).
\end{enumerate}
\end{thm}
\begin{thm}
\label{thm:condition_for_finite_generation_without_cycle}
Let $r \in S$.
Suppose that $(\Gamma^{\mathrm{odd}})_{\sim r}$ does not contain a nontrivial cycle.
Then $W^{\perp r}$ is finitely generated if and only if all the following four conditions are satisfied (see above for the notations):
\begin{enumerate}
\item \label{item:thm_condition_for_finite_generation_without_cycle__finite_order_almost_everywhere}
There exist only finitely many pairs $(x,y)$ of elements $x,y \in \mathrm{Odd}(r)$ for which $m(x,y)$ is an even integer.
\item \label{item:thm_condition_for_finite_generation_without_cycle__E_finite}
The set $\partial\mathrm{EOdd}(r)$ is finite.
\item \label{item:thm_condition_for_finite_generation_without_cycle__2_or_infinite_almost_everywhere}
For each $x \in \partial\mathrm{EOdd}(r)$, there exist only finitely many elements $y \in \mathrm{Odd}(r)$ for which $m(x,y) \not\in \{2,\infty\}$.
\item \label{item:thm_condition_for_finite_generation_without_cycle__finitely_many_components}
For each $x \in \partial\mathrm{EOdd}(r)$, the graph $\Gamma_{\mathrm{Odd}(r)^{\perp x}}^{\mathrm{odd}}$ consists of finitely many connected components.
\end{enumerate}
\end{thm}

Here we present some corollaries of the above main theorems applied to some subclasses of Coxeter systems.
The first case is that the Coxeter group $W$ itself is finitely generated:
\begin{cor}
\label{cor:condition_for_finite_generation__W_finitely_generated}
Suppose that $W$ is finitely generated.
Let $r \in S$.
Then $W^{\perp r}$ is finitely generated if and only if either $(\Gamma^{\mathrm{odd}})_{\sim r}$ does not contain a nontrivial cycle, or $(\Gamma^{\mathrm{odd}})_{\sim r}$ contains a nontrivial cycle and Conditions \ref{item:thm_condition_for_finite_generation_with_cycle__core}, \ref{item:thm_condition_for_finite_generation_with_cycle__2_or_infinite} and \ref{item:thm_condition_for_finite_generation_with_cycle__non-adjacent_vertices} in Theorem \ref{thm:condition_for_finite_generation_with_cycle} are satisfied.
\end{cor}
\begin{proof}
By the hypothesis $|S| < \infty$ here, Condition \ref{item:thm_condition_for_finite_generation_with_cycle__E_finite} in Theorem \ref{thm:condition_for_finite_generation_with_cycle} and all of Conditions \ref{item:thm_condition_for_finite_generation_without_cycle__finite_order_almost_everywhere}--\ref{item:thm_condition_for_finite_generation_without_cycle__finitely_many_components} in Theorem \ref{thm:condition_for_finite_generation_without_cycle} are automatically satisfied.
\end{proof}
In particular, the reflection subgroup $W^{\perp r}$ is always finitely generated if $W$ is an affine Coxeter group.
This is a special case of a general theorem of Dyer \cite[Theorem 5.1]{Dye} which proves that any reflection subgroup of an affine Coxeter group is finitely generated.

A Coxeter system $(W,S)$ is called \emph{2-spherical} if we have $m(s,t) < \infty$ for every $s,t \in S$.
For this case, we have the following:
\begin{cor}
\label{cor:condition_for_finite_generation__2-spherical}
Suppose that $(W,S)$ is 2-spherical and irreducible.
Let $r \in S$.
Then $W^{\perp r}$ is finitely generated if and only if one of the following three conditions is satisfied:
\begin{enumerate}
\item \label{item:cor_condition_for_finite_generation__2-spherical__odd_everywhere}
$m(x,y)$ is an odd integer for every $x,y \in S$;
\item \label{item:cor_condition_for_finite_generation__2-spherical__A_tilde}
$\Gamma$ is of type $\widetilde{A}_{n-1}$ with $4 \leq n$ (in particular, $|S| < \infty$);
\item \label{item:cor_condition_for_finite_generation__2-spherical__acyclic}
$|S| < \infty$ and $(\Gamma^{\mathrm{odd}})_{\sim r}$ does not contain a nontrivial cycle.
\end{enumerate}
\end{cor}
\begin{proof}
First, note that we have $S = \mathrm{EOdd}(r)$ in this case.
The \lq\lq if'' part of the claim is easily proven by verifying the conditions in Theorems \ref{thm:condition_for_finite_generation_with_cycle} and \ref{thm:condition_for_finite_generation_without_cycle} (see also Corollary \ref{cor:condition_for_finite_generation__W_finitely_generated} for Case \ref{item:cor_condition_for_finite_generation__2-spherical__acyclic}).
From now, we consider the \lq\lq only if'' part.

First we consider the case that $(\Gamma^{\mathrm{odd}})_{\sim r}$ does not contain a nontrivial cycle as in Theorem \ref{thm:condition_for_finite_generation_without_cycle}.
We show that now $|S| < \infty$, which leads to Condition \ref{item:cor_condition_for_finite_generation__2-spherical__acyclic}.
Since $S = \mathrm{EOdd}(r)$ as above and $\partial\mathrm{EOdd}(r)$ is finite by Property \ref{item:thm_condition_for_finite_generation_without_cycle__E_finite} in Theorem \ref{thm:condition_for_finite_generation_without_cycle}, it suffices to show that $|\mathrm{Odd}(r)| < \infty$.
Let $I_{\mathrm{odd}}$ (respectively, $I_{\mathrm{even}}$) be the subsets of $\mathrm{Odd}(r) \smallsetminus \{r\}$ consisting of the elements $x$ for which $m(r,x)$ is an odd (respectively, even) integer.
We have $\mathrm{Odd}(r) = I_{\mathrm{odd}} \cup I_{\mathrm{even}} \cup \{r\}$ since $(W,S)$ is 2-spherical.
Now Property \ref{item:thm_condition_for_finite_generation_without_cycle__finite_order_almost_everywhere} implies that $|I_{\mathrm{even}}| < \infty$.
On the other hand, since $(W,S)$ is 2-spherical and $(\Gamma^{\mathrm{odd}})_{\sim r}$ is assumed to have no nontrivial cycles, $m(x,y)$ is an even integer for every distinct $x,y \in I_{\mathrm{odd}}$.
Therefore, Property \ref{item:thm_condition_for_finite_generation_without_cycle__finite_order_almost_everywhere} implies that $|I_{\mathrm{odd}}| < \infty$.
Hence we have $|\mathrm{Odd}(r)| < \infty$, as desired.

Secondly, we consider the case that $(\Gamma^{\mathrm{odd}})_{\sim r}$ contains a nontrivial cycle as in Theorem \ref{thm:condition_for_finite_generation_with_cycle}.
Since $(W,S)$ is 2-spherical, Property \ref{item:thm_condition_for_finite_generation_with_cycle__2_or_infinite} in Theorem \ref{thm:condition_for_finite_generation_with_cycle} implies that the sets $\mathrm{Odd}(r)$ and $\partial\mathrm{EOdd}(r)$ are not adjacent in $\Gamma$; therefore we have $\partial\mathrm{EOdd}(r) = \emptyset$ and $S = \mathrm{Odd}(r)$ since $S = \mathrm{EOdd}(r)$ is irreducible.
Now, since $(W,S)$ is 2-spherical, Property \ref{item:thm_condition_for_finite_generation_with_cycle__core_not_even} leads to Condition \ref{item:cor_condition_for_finite_generation__2-spherical__odd_everywhere} in the statement, while Property \ref{item:thm_condition_for_finite_generation_with_cycle__core_bipyramid} cannot be satisfied.
Finally, we consider the case that Property \ref{item:thm_condition_for_finite_generation_with_cycle__core_A_tilde} holds.
Now if $\mathrm{Odd}(r) \smallsetminus K \neq \emptyset$, then there exists an element $x \in \mathrm{Odd}(r) \smallsetminus K$ which is adjacent to some $y \in K$ in the graph $\Gamma^{\mathrm{odd}}$.
Moreover, by the type of $K$ specified in Property \ref{item:thm_condition_for_finite_generation_with_cycle__core_A_tilde}, there exists an element $z \in K \smallsetminus \{y\}$ for which $m(y,z)$ is an odd integer.
Now $m(x,z)$ cannot be an odd integer; otherwise, a cycle $(x,z,y,x)$ in $(\Gamma^{\mathrm{odd}})_{\sim r}$ is not contained in the core $K$ of $(\Gamma^{\mathrm{odd}})_{\sim r}$, a contradiction.
Therefore $m(x,z)$ is an even integer since $(W,S)$ is 2-spherical, contradicting Property \ref{item:thm_condition_for_finite_generation_with_cycle__core_A_tilde}.
Hence we have $K = \mathrm{Odd}(r) = S$, which leads to Condition \ref{item:cor_condition_for_finite_generation__2-spherical__A_tilde}.
This completes the proof of Corollary \ref{cor:condition_for_finite_generation__2-spherical}.
\end{proof}

A Coxeter system $(W,S)$ is called \emph{even} if $m(s,t)$ is an even integer or $\infty$ for every distinct $s,t \in S$.
For this case, we have the following:
\begin{cor}
\label{cor:condition_for_finite_generation__even}
Suppose that $(W,S)$ is even.
Let $r \in S$.
Then $W^{\perp r}$ is finitely generated if and only if there exist only finitely many elements $x \in S$ for which $m(r,x) < \infty$.
\end{cor}
\begin{proof}
Note that $\mathrm{Odd}(r) = \{r\}$ since $(W,S)$ is even.
Therefore, it is the situation in Theorem \ref{thm:condition_for_finite_generation_without_cycle}, and Conditions \ref{item:thm_condition_for_finite_generation_without_cycle__finite_order_almost_everywhere}, \ref{item:thm_condition_for_finite_generation_without_cycle__2_or_infinite_almost_everywhere} and \ref{item:thm_condition_for_finite_generation_without_cycle__finitely_many_components} are automatically satisfied.
This implies that $W^{\perp r}$ is finitely generated if and only if $\partial\mathrm{EOdd}(r)$ is a finite set.
Now we have $\partial\mathrm{EOdd}(r) = \{x \in S \smallsetminus \{r\} \mid m(r,x) < \infty\}$ since $(W,S)$ is even, therefore the claim follows.
\end{proof}

A Coxeter system $(W,S)$ is called \emph{skew-angled} if we have $m(s,t) \neq 2$ for every $s,t \in S$.
For this case, we have the following:
\begin{cor}
\label{cor:condition_for_finite_generation__skew-angled}
Suppose that $(W,S)$ is skew-angled.
Let $r \in S$, and suppose further that $(\Gamma^{\mathrm{odd}})_{\sim r}$ contains a nontrivial cycle.
Then $W^{\perp r}$ is finitely generated if and only if $m(x,y)$ is an odd integer or $\infty$ for every $x,y \in \mathrm{Odd}(r)$, and $m(x,y) = \infty$ for every $x \in \mathrm{Odd}(r)$ and $y \in S \smallsetminus \mathrm{Odd}(r)$ (i.e., $W$ is the free product of $W_{\mathrm{Odd}(r)}$ and $W_{S \smallsetminus \mathrm{Odd}(r)}$).
\end{cor}
\begin{proof}
For the \lq\lq if'' part, we have $\partial\mathrm{EOdd}(r) = \emptyset$ by the second condition, while Condition \ref{item:thm_condition_for_finite_generation_with_cycle__core_not_even} in Theorem \ref{thm:condition_for_finite_generation_with_cycle} follows from the first condition.
Hence the \lq\lq if'' part holds.
For the \lq\lq only if'' part, Properties \ref{item:thm_condition_for_finite_generation_with_cycle__core_A_tilde} and \ref{item:thm_condition_for_finite_generation_with_cycle__core_bipyramid} in Theorem \ref{thm:condition_for_finite_generation_with_cycle} cannot hold since $(W,S)$ is skew-angled, therefore Property \ref{item:thm_condition_for_finite_generation_with_cycle__core_not_even} (hence the first condition in the statement) holds.
Moreover, Property \ref{item:thm_condition_for_finite_generation_with_cycle__2_or_infinite} now implies that $m(x,y) = \infty$ for every $x \in \mathrm{Odd}(r)$ and $y \in \partial\mathrm{EOdd}(r)$, therefore we have $\partial\mathrm{EOdd}(r) = \emptyset$ and the second condition in the statement holds.
Hence the claim follows.
\end{proof}

Finally, we consider the following condition for $(W,S)$; every subset $I \subseteq S$ with at least $3$ elements generates an infinite subgroup.
(This is equivalent to the property that the Davis-Vinberg complex associated to $(W,S)$ has dimension at most two.)
For this case, we have the following:
\begin{cor}
\label{cor:condition_for_finite_generation__2-dimensional}
Suppose that every subset $I \subseteq S$ with $|I| \geq 3$ generates an infinite subgroup of $W$.
Let $r \in S$, and suppose further that $(\Gamma^{\mathrm{odd}})_{\sim r}$ contains a nontrivial cycle.
Then $W^{\perp r}$ is finitely generated if and only if $m(x,y)$ is an odd integer or $\infty$ for every $x,y \in \mathrm{Odd}(r)$, and $m(x,y) = \infty$ for every $x \in \mathrm{Odd}(r)$ and $y \in S \smallsetminus \mathrm{Odd}(r)$ (i.e., $W$ is the free product of $W_{\mathrm{Odd}(r)}$ and $W_{S \smallsetminus \mathrm{Odd}(r)}$).
\end{cor}
\begin{proof}
The proof of \lq\lq if'' part is the same as Corollary \ref{cor:condition_for_finite_generation__skew-angled}.
For the \lq\lq only if'' part, the hypothesis implies that there do not exist three elements $z_1,z_2,z_3 \in \mathrm{Odd}(r)$ satisfying that $m(z_1,z_2) = m(z_2,z_3) = 3$ and $m(z_1,z_3) = 2$ (i.e., generating a parabolic subgroup of type $A_3$).
This implies that Properties \ref{item:thm_condition_for_finite_generation_with_cycle__core_A_tilde} and \ref{item:thm_condition_for_finite_generation_with_cycle__core_bipyramid} in Theorem \ref{thm:condition_for_finite_generation_with_cycle} cannot hold, therefore Property \ref{item:thm_condition_for_finite_generation_with_cycle__core_not_even} (hence the first condition in the statement) holds.
Now our remaining task is to show that $\partial\mathrm{EOdd}(r) = \emptyset$.
Assume contrary that $\partial\mathrm{EOdd}(r)$ contains an element $x$.
Then Property \ref{item:thm_condition_for_finite_generation_with_cycle__non-adjacent_vertices} in Theorem \ref{thm:condition_for_finite_generation_with_cycle} implies that $\mathrm{Odd}(r)^{\perp x}$ is a pre-core of $(\Gamma^{\mathrm{odd}})_{\sim r}$; in particular, there exist two distinct elements $y,z \in \mathrm{Odd}(r)^{\perp x}$ for which $m(y,z)$ is an odd integer.
However, now the subset $\{x,y,z\}$ of $S$ generates a finite parabolic subgroup of type $A_1 \times I_2(m(y,z))$ (where we put $I_2(3) = A_2$), contradicting the assumption on $(W,S)$.
Hence we have $\partial\mathrm{EOdd}(r) = \emptyset$, concluding the proof.
\end{proof}

\subsection{Proof of main theorems}
\label{sec:main_theorem_proof}

In this subsection, we prove Theorem \ref{thm:condition_for_finite_generation_with_cycle} and Theorem \ref{thm:condition_for_finite_generation_without_cycle}.
Note that by Lemma \ref{lem:finitely_generated_iff_finite_rank}, for $r \in S$, $W^{\perp r}$ is finitely generated if and only if $R^{r}$ is a finite set.

\begin{figure}[hbt]
\centering
\begin{picture}(100,70)(0,-70)
\put(20,-35){\circle{6}}\put(50,-5){\circle{6}}\put(50,-65){\circle{6}}\put(80,-35){\circle{6}}
\put(22,-33){\line(1,1){26}}\put(22,-37){\line(1,-1){26}}\put(78,-33){\line(-1,1){26}}\put(78,-37){\line(-1,-1){26}}\put(50,-8){\line(0,-1){54}}
\put(53,-37){$\infty$}
\end{picture}
\caption{Coxeter graph for the latter case of Lemma \ref{lem:shape_of_cycle}}
\label{fig:diamond}
\end{figure}
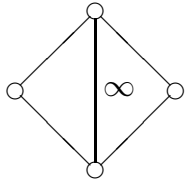

We start with a few lemmas on some substructures of $W^{\perp r}$:
\begin{lem}
\label{lem:no_realizable}
Let $r \in I \subseteq S$.
Suppose that $(\Gamma_I^{\mathrm{odd}})_{\sim r}$ contains a nontrivial cycle, and an element $\xi \in \mathcal{E}^I(r)$ satisfies the following condition:
\begin{equation}
\begin{split}
\label{eq:lem_no_realizable_condition}
&\mbox{No element } 1_{\xi^{\circ}} \neq c \in \pi_1(\Gamma_I^{\mathrm{odd}};\xi^{\circ}) \mbox{ is realized by} \\
&\mbox{a sequence of moves in } I \mbox{ from } \xi \mbox{ to } \xi \enspace.
\end{split}
\end{equation}
Then all the elements $s_r^I(c;\xi) \in R^{I,r}$ with $c \in \pi_1(\Gamma_I^{\mathrm{odd}};r,\xi^{\circ})$ are different from each other, hence we have $|R^{I,r}| = \infty$.
\end{lem}
\begin{proof}
By virtue of the condition (\ref{eq:lem_no_realizable_condition}), the first part of the claim follows from Theorem \ref{thm:presentation_of_orthogonal_reflection_subgroup}.
Moreover, by the hypothesis, the free group $\pi_1(\Gamma_I^{\mathrm{odd}};r)$ is nontrivial, therefore we have $|\pi_1(\Gamma_I^{\mathrm{odd}};r,\xi^{\circ})| = |\pi_1(\Gamma_I^{\mathrm{odd}};r)| = \infty$.
Hence the second part of the claim follows from the first part.
\end{proof}
\begin{lem}
\label{lem:adjacent_closed_path}
Let $r \in I \subseteq S$ and $\xi \in \mathcal{E}^I(r)$ with $\xi^{\circ} = r$.
Suppose that $P = (x_n,\dots,x_1,x_0)$ is a nontrivial reduced (but not necessarily cyclically reduced) closed path in $\Gamma_I^{\mathrm{odd}}$ with $x_0 = x_n = r$, $\xi_{\dagger} \not\in V(P)$ and $\xi_{\dagger}$ is adjacent to $V(P)$ in $\Gamma$.
Then all the roots $\gamma_r^I([P]^k;\xi) \in \Pi^{I,r}$ with $k \in \mathbb{Z}$ are different from each other, hence we have $|R^{I,r}| = \infty$.
\end{lem}
\begin{proof}
First, note that $\xi_{\dagger} \in \mathrm{supp}(\gamma(\xi)) \subseteq \{\xi^{\circ},\xi_{\dagger}\}$ by the property of $\gamma(\xi)$.
Secondly, since the closed path $P$ is reduced, there exists an index $\ell \geq 1$ with $x_{\ell} \neq x_{n-\ell}$.
Take the minimal $\ell$ with this property.
Then $P$ can be decomposed as $P = Q^{-1}CQ$ where $Q = (x_{\ell-1},\dots,x_1,x_0)$ is a reduced path and $C = (x_{n-\ell+1},n_{n-\ell},\dots,x_{\ell},x_{\ell-1})$ is a nontrivial, cyclically reduced closed path.
Now for each integer $k \neq 0$, we have $[P]^k = [Q^{-1}C^kQ]$, $V(Q^{-1}C^kQ) = V(P)$ and $Q^{-1}C^kQ$ is a nontrivial reduced closed path.
Therefore we have $\mathrm{supp}(\pi^I([P]^k)) = V(P)$ by Corollary \ref{cor:support_of_reduced_path}.
Since $\xi_{\dagger} \not\in V(P)$ and $\xi_{\dagger}$ is adjacent to $V(P)$ in $\Gamma$ by the hypothesis, we have $\pi^I([P]^k) \cdot \alpha_{\xi_{\dagger}} \neq \alpha_{\xi_{\dagger}}$ by Lemma \ref{lem:wdoesnotfixgamma}, while we have $\pi^I([P]^k) \cdot \alpha_{\xi^{\circ}} = \alpha_{\xi^{\circ}}$ since $\pi^I([P]^k) \in Y_{\xi^{\circ},\xi^{\circ}}^I$.
This implies that $\pi^I([P]^k) \cdot \gamma(\xi) \neq \gamma(\xi)$, therefore we have $\gamma_r^I([P]^i;\xi) \neq \gamma_r^I([P]^j;\xi)$ for any distinct $i,j \in \mathbb{Z}$.
Hence the claim holds.
\end{proof}
\begin{lem}
\label{lem:shape_of_cycle}
Let $r \in I \subseteq S$ and $P = (x_n,\dots,x_1,x_0)$ be a nontrivial full cycle in $(\Gamma_I^{\mathrm{odd}})_{\sim r}$.
Suppose that $m(x_i,x_j)$ is an even integer for some indices $i,j$ (hence $n \geq 4$).
If $|R^{I,r}| < \infty$, then we have either $V(P)$ is of type $\widetilde{A}_{n-1}$, or $n = 4$, $m(x_0,x_1) = m(x_1,x_2) = m(x_2,x_3) = m(x_3,x_0) = 3$, $m(x_0,x_2) = 2$ and $m(x_1,x_3) = \infty$ (see Figure \ref{fig:diamond}).
\end{lem}
\begin{proof}
Take indices $i,j$ for which $m(x_i,x_j)$ is an even integer.
By Remark \ref{rem:restriction_to_smaller_domain} and Remark \ref{rem:isomorphism_by_conjugation}, we may assume without loss of generality that $r = x_i$ and $I = V(P)$.
Now $\pi_1(\Gamma_I^{\mathrm{odd}};r)$ is the infinite cyclic group generated by $[P^{\to i}]$.
By the hypothesis $|R^{I,r}| < \infty$, Theorem \ref{thm:presentation_of_orthogonal_reflection_subgroup} implies that at least one element $[P^{\to}]^{\ell}$ with $0 \neq \ell \in \mathbb{Z}$ can be realized by a sequence of moves in $I$ from $(x_i,x_j) \in \mathcal{E}^I(r)$ to $(x_i,x_j)$.
By virtue of Example \ref{ex:cycle_odd_length} and Example \ref{ex:cycle_even_length}, this is possible only when either $V(P)$ is of type $\widetilde{A}_{n-1}$, or we have $n = 4$, $m(x_0,x_1) = m(x_1,x_2) = m(x_2,x_3) = m(x_3,x_0) = 3$, $m(x_0,x_2) = 2$ and $m(x_1,x_3)$ is either $\infty$ or an even integer not equal to $2$.
Moreover, in the latter case with $2 < m(x_1,x_3) < \infty$, we have $(x_1,x_3) \in \mathcal{E}^I(r)$, and the same argument implies that no nontrivial element of the infinite group $\pi_1(\Gamma_I^{\mathrm{odd}};x_1)$ can be realized by a sequence of moves in $I$ from $(x_1,x_3)$ to $(x_1,x_3)$.
Since $|R^{I,r}| < \infty$, this contradicts Theorem \ref{thm:presentation_of_orthogonal_reflection_subgroup}.
Hence the claim holds.
\end{proof}

The next proposition is a key observation in our argument:
\begin{prop}
\label{prop:cycle_contains_even_edges}
Let $r \in I \subseteq S$ and $P = (x_n,\dots,x_1,x_0)$ be a nontrivial cycle in $(\Gamma_I^{\mathrm{odd}})_{\sim r}$.
If $|R^{I,r}| < \infty$, then for each $\xi \in \mathcal{E}^I(r)$ with $\xi_{\dagger} \in \mathrm{Odd}_I(r)$, we have $m(\xi^{\circ},\xi_{\dagger}) = 2$ and $\xi^{\circ},\xi_{\dagger} \in V(P)$.
\end{prop}
\begin{proof}
We prove the contraposition that $|R^{I,r}| = \infty$ if there exists an element $\xi \in \mathcal{E}^I(r)$ with $\xi_{\dagger} \in \mathrm{Odd}_I(r)$ satisfying either $m(\xi^{\circ},\xi_{\dagger}) > 2$ or $\{\xi^{\circ},\xi_{\dagger}\} \not\subseteq V(P)$.
In the case $m(\xi^{\circ},\xi_{\dagger}) > 2$, $\xi$ admits no moves in $I$, therefore the claim follows from Lemma \ref{lem:no_realizable}.
From now, we consider the case that $\{\xi^{\circ},\xi_{\dagger}\} \not\subseteq V(P)$.

We use induction on the distance $d(\xi^{\circ},V(P))$ in $\Gamma_I^{\mathrm{odd}}$ from $\xi^{\circ}$ to $V(P)$; and for each case, we use induction on the distance $d(\xi_{\dagger},V(P))$ in $\Gamma_I^{\mathrm{odd}}$ from $\xi_{\dagger}$ to $V(P)$.
We may assume without loss of generality that $d(\xi^{\circ},V(P)) \leq d(\xi_{\dagger},V(P))$, since otherwise we can consider another element $(\xi_{\dagger},\xi^{\circ}) \in \mathcal{E}^I(r)$ instead of $\xi$.
Let $Q$ be one of the shortest paths in $\Gamma_I^{\mathrm{odd}}$ from $\xi^{\circ}$ to a vertex, say $x_k$, in $V(P)$ (hence $V(Q) \cap V(P) = \{x_k\}$).
By the above-mentioned relation of distances, $\xi_{\dagger}$ is not contained in the reduced closed path $Q^{-1}P^{\to k}Q$ in $\Gamma_I^{\mathrm{odd}}$.
If $\xi_{\dagger}$ is adjacent to $V(Q^{-1}P^{\to k}Q) = V(P) \cup V(Q)$ in $\Gamma$, then we have $|R^{I,\xi^{\circ}}| = \infty$ by Lemma \ref{lem:adjacent_closed_path}, therefore the claim follows from Remark \ref{rem:isomorphism_by_conjugation}.
From now, we consider the other case that $\xi_{\dagger}$ is not adjacent to $V(P) \cup V(Q)$ in $\Gamma$.

If $d(\xi^{\circ},V(P)) > 0$, and $y$ is the vertex in $Q$ next to $\xi^{\circ}$, then (since $m(y,\xi_{\dagger}) = 2$ as above) the claim follows from the induction hypothesis applied to $(y,\xi_{\dagger}) \in \mathcal{E}^I(r)$.
From now, we consider the remaining case that $d(\xi^{\circ},V(P)) = 0$, i.e., $\xi^{\circ} \in V(P)$.

Let $Q'$ be one of the shortest paths in $\Gamma_I^{\mathrm{odd}}$ from $\xi_{\dagger}$ to a vertex in $V(P)$.
Since $\xi_{\dagger}$ is not adjacent to $V(P)$ in $\Gamma$, $Q'$ has at least two edges, and we may assume without loss of generality that $\xi^{\circ}$ is not the end vertex of $Q'$.
Let $y$ be the vertex in $Q'$ next to $\xi_{\dagger}$ (hence $y \not\in V(P)$).
Now if $m(\xi^{\circ},y)$ is an odd integer, then the concatenation of the edge from $y$ to $\xi^{\circ}$, a subpath in $P$ from $\xi^{\circ}$ to the end vertex of $Q'$, and the subpath in $Q'{}^{-1}$ from the end vertex of $Q'$ to $y$ forms a nontrivial cycle in $\Gamma_I^{\mathrm{odd}}$ containing $\xi^{\circ}$ and having smaller distance to $\xi_{\dagger}$ than $P$; hence the claim follows by induction hypothesis applied to the cycle and the element $\xi$.
On the other hand, if $m(\xi^{\circ},y)$ is an even integer, then the claim follows from the induction hypothesis applied to $(\xi^{\circ},y) \in \mathcal{E}^I(r)$.
Finally, suppose that $m(\xi^{\circ},y) = \infty$.
Put $J = V(P) \cup V(Q')$.
In this case, $\xi_{\dagger}$ is not adjacent in $\Gamma_J^{\mathrm{odd}}$ to any vertex other than $y$, therefore $(\xi_{\dagger},\xi^{\circ}) \in \mathcal{E}^J(\xi_{\dagger})$ admits no moves in $J$.
Now Lemma \ref{lem:no_realizable} implies that $|R^{J,\xi_{\dagger}}| = \infty$, therefore we have $|R^{I,r}| = \infty$ by Remark \ref{rem:restriction_to_smaller_domain} and Remark \ref{rem:isomorphism_by_conjugation}, as desired.
Hence the proof of Proposition \ref{prop:cycle_contains_even_edges} is concluded.
\end{proof}

We introduce some auxiliary terminology.
Let $r \in S$.
For an element $s_r(c;\xi)$ of $R^{r}$, we say that it is an \emph{inner generator} (respectively, \emph{outer generator}) if $\xi_{\dagger} \in \mathrm{Odd}(r)$ (respectively, $\xi_{\dagger} \not\in \mathrm{Odd}(r)$).
By Theorem \ref{thm:presentation_of_orthogonal_reflection_subgroup} and the definitions of S-moves and E-moves, the properties of being an inner generator and of being an outer generator are independent of the expression $s_r(c;\xi)$ of an element of $R^{r}$.
Note that we have $\xi_{\dagger} \in \partial\mathrm{EOdd}(r)$ for each outer generator $s_r(c;\xi)$.
Note also that for any outer generator $s_r(c;\xi)$, the element $\xi \in \mathcal{E}(r)$ admits no E-moves by the definition.

First, we determine the condition for the number of outer generators to be finite:
\begin{prop}
\label{prop:when_outer_generators_are_finite}
Let $r \in S$.
If $(\Gamma^{\mathrm{odd}})_{\sim r}$ contains a nontrivial cycle, then the number of outer generators is finite if and only if Conditions \ref{item:thm_condition_for_finite_generation_with_cycle__E_finite}--\ref{item:thm_condition_for_finite_generation_with_cycle__non-adjacent_vertices} in Theorem \ref{thm:condition_for_finite_generation_with_cycle} are satisfied.
On the other hand, if $(\Gamma^{\mathrm{odd}})_{\sim r}$ does not contain a nontrivial cycle, then the number of outer generators is finite if and only if Conditions \ref{item:thm_condition_for_finite_generation_without_cycle__E_finite}--\ref{item:thm_condition_for_finite_generation_without_cycle__finitely_many_components} are satisfied.
\end{prop}
\begin{proof}
First note that, by the definitions of S-moves and E-moves, Theorem \ref{thm:presentation_of_orthogonal_reflection_subgroup} implies that two outer generators $s_r(c;\xi)$ and $s_r(d;\zeta)$ cannot coincide if $\xi_{\dagger} \neq \zeta_{\dagger}$.
Hence the condition $|\partial\mathrm{EOdd}(r)| < \infty$ is necessary for the number of outer generators to be finite.
From now, we suppose that $|\partial\mathrm{EOdd}(r)| < \infty$.

We consider the case that $(\Gamma^{\mathrm{odd}})_{\sim r}$ contains a nontrivial cycle.
First we prove the \lq\lq if'' part in the statement.
Let $s_r(c;\xi)$ and $s_r(d;\zeta)$ be outer generators with $\xi_{\dagger} = \zeta_{\dagger}$.
By Condition \ref{item:thm_condition_for_finite_generation_without_cycle__2_or_infinite_almost_everywhere} in Theorem \ref{thm:condition_for_finite_generation_with_cycle}, we have $\xi^{\circ},\zeta^{\circ} \in K' := \mathrm{Odd}(r)^{\perp \xi_{\dagger}}$.
Now by Condition \ref{item:thm_condition_for_finite_generation_with_cycle__non-adjacent_vertices} in Theorem \ref{thm:condition_for_finite_generation_with_cycle}, the set $K'$ is a pre-core of $(\Gamma^{\mathrm{odd}})_{\sim r}$, therefore we have $c^{-1}d \in \pi_1(\Gamma^{\mathrm{odd}};\xi^{\circ},\zeta^{\circ}) = \pi_1(\Gamma_{K'}^{\mathrm{odd}};\xi^{\circ},\zeta^{\circ})$ by Lemma \ref{lem:fundamental_group_for_subgraph}.
This implies that the equivalence class $c^{-1}d$ is represented by a path in $\Gamma_{K'}^{\mathrm{odd}}$ from $\zeta^{\circ}$ to $\xi^{\circ}$, which (by the choice of $K'$ and the definition of S-moves) can be realized by a sequence of S-moves from $\zeta$ to $\xi$.
Hence we have $s_r(c;\xi) = s_r(d;\zeta)$ by Theorem \ref{thm:presentation_of_orthogonal_reflection_subgroup}.
Therefore, the number of outer generators is equal to $|\partial\mathrm{EOdd}(r)| < \infty$.

Secondly, we prove the \lq\lq only if'' part.
If $m(x,y) \not\in \{2,\infty\}$ for some $x \in \mathrm{Odd}(r)$ and $y \in \partial\mathrm{EOdd}(r)$, then $m(x,y)$ is an even integer different from $2$.
Now we have $(x,y) \in \mathcal{E}(r)$ and $(x,y)$ admits no moves, therefore $|R^{r}| = \infty$ by Lemma \ref{lem:no_realizable}.
This is a contradiction, therefore Condition \ref{item:thm_condition_for_finite_generation_with_cycle__2_or_infinite} in Theorem \ref{thm:condition_for_finite_generation_with_cycle} is satisfied.
On the other hand, assume contrary that $\mathrm{Odd}(r)^{\perp x}$ is not a pre-core of $(\Gamma^{\mathrm{odd}})_{\sim r}$ for some $x \in \partial\mathrm{EOdd}(r)$.
Then there exist an element $y \in \mathrm{Odd}(r)^{\perp x}$ and a nontrivial cycle $P$ in $(\Gamma^{\mathrm{odd}})_{\sim r}$ satisfying that $V(P)$ is not contained in $(\Gamma_{\mathrm{Odd}(r)^{\perp x}}^{\mathrm{odd}})_{\sim y}$ (note that $\mathrm{Odd}(r)^{\perp x} \neq \emptyset$ by definition of $\partial\mathrm{EOdd}(r)$ and the above result).
Take a shortest path $Q$ in $\Gamma^{\mathrm{odd}}$ from $y$ to $V(P)$, and take a cyclic shift $P'$ of $P$ starting at the end vertex of $Q$.
Then $Q^{-1}P'Q$ is a nontrivial reduced closed path starting at $y$ and $V(Q^{-1}P'Q) \not\subseteq \mathrm{Odd}(r)^{\perp x}$.
Now Lemma \ref{lem:adjacent_closed_path} implies that $|R^{y}| = \infty$, therefore $|R^{r}| = \infty$ by Remark \ref{rem:isomorphism_by_conjugation}, a contradiction.
Hence Condition \ref{item:thm_condition_for_finite_generation_with_cycle__non-adjacent_vertices} in Theorem \ref{thm:condition_for_finite_generation_with_cycle} is satisfied, concluding the proof for the case that $(\Gamma^{\mathrm{odd}})_{\sim r}$ contains a nontrivial cycle.

From now, we consider the other case that $(\Gamma^{\mathrm{odd}})_{\sim r}$ does not contain a nontrivial cycle.
Let $s_r(c;\xi)$ and $s_r(d;\zeta)$ be outer generators with $\xi_{\dagger} = \zeta_{\dagger}$.
In this case, Theorem \ref{thm:presentation_of_orthogonal_reflection_subgroup} implies that we have $s_r(c;\xi) = s_r(d;\zeta)$ if and only if there exists a sequence of S-moves from $\xi$ to $\zeta$.
Now, since $\xi$ admits no S-moves when $m(\xi^{\circ},\xi_{\dagger}) > 2$, it follows that the number of outer generators $s_r(c;\xi)$ with $m(\xi^{\circ},\xi_{\dagger}) > 2$ is finite if and only if Condition \ref{item:thm_condition_for_finite_generation_without_cycle__2_or_infinite_almost_everywhere} in Theorem \ref{thm:condition_for_finite_generation_without_cycle} is satisfied.
On the other hand, by the definition of S-moves, for the case that $m(\xi^{\circ},\xi_{\dagger}) = m(\zeta^{\circ},\zeta_{\dagger}) = 2$, we have $s_r(c;\xi) = s_r(d;\zeta)$ if and only if $\xi^{\circ}$ and $\zeta^{\circ}$ belong to the same connected component of $\Gamma_{\mathrm{Odd}(r)^{\perp \xi_{\dagger}}}^{\mathrm{odd}}$.
Therefore, the number of outer generators $s_r(c;\xi)$ with $m(\xi^{\circ},\xi_{\dagger}) = 2$ is finite if and only if Condition \ref{item:thm_condition_for_finite_generation_without_cycle__finitely_many_components} in Theorem \ref{thm:condition_for_finite_generation_without_cycle} is satisfied.
Hence the proof of Proposition \ref{prop:when_outer_generators_are_finite} is concluded.
\end{proof}

Secondly, we determine the condition for the number of inner generators to be finite.
The next proposition deals with the case that $(\Gamma^{\mathrm{odd}})_{\sim r}$ does not contain a nontrivial cycle:
\begin{prop}
\label{prop:when_inner_generators_are_finite__acyclic}
Let $r \in S$, and suppose that $(\Gamma^{\mathrm{odd}})_{\sim r}$ does not contain a nontrivial cycle.
Then the number of inner generators is finite if and only if Condition \ref{item:thm_condition_for_finite_generation_without_cycle__finite_order_almost_everywhere} in Theorem \ref{thm:condition_for_finite_generation_without_cycle} is satisfied.
\end{prop}
\begin{proof}
For the \lq\lq if'' part, note that $|\mathcal{E}^{\mathrm{Odd}(r)}(r)| = |\mathcal{R}^{\mathrm{Odd}(r)}(r)|$ by the hypothesis on the shape of $(\Gamma^{\mathrm{odd}})_{\sim r}$, while $\mathcal{E}^{\mathrm{Odd}(r)}(r)$ is finite by Condition \ref{item:thm_condition_for_finite_generation_without_cycle__finite_order_almost_everywhere}.
This implies that $\mathcal{R}^{\mathrm{Odd}(r)}(r)$ is finite, so is $R^{\mathrm{Odd}(r),r}$ (see Theorem \ref{thm:generator_of_orthogonal_reflections}), as desired.

Secondly, for the \lq\lq only if'' part, let $s_r(c_1;\xi_1),\dots,s_r(c_k;\xi_k)$ be all of the inner generators.
Then by Theorem \ref{thm:presentation_of_orthogonal_reflection_subgroup}, for each pair $x,y$ of elements of $\mathrm{Odd}(r)$ for which $m(x,y)$ is an even integer, there exists a sequence of moves from $(x,y) \in \mathcal{E}(r)$ to one of the $\xi_i$.
Now for each $1 \leq i \leq k$, let $P_i$ denote the unique reduced path in $\Gamma^{\mathrm{odd}}$ from $(\xi_i)_{\dagger}$ to $(\xi_i)^{\circ}$ (recall the hypothesis on the shape of $(\Gamma^{\mathrm{odd}})_{\sim r}$).
Let $z_i$ denote the vertex in $P_i$ next to $(\xi_i)_{\dagger}$, and let $z'_i$ denote the vertex in $P_i$ next to $z_i$.
Then by the shape of $(\Gamma^{\mathrm{odd}})_{\sim r}$ and the definitions of S-moves and E-moves, the subset of $\mathcal{E}(r)$ consisting of elements of the form $(x,(\xi_i)_{\dagger})$ with $x \in \mathrm{Odd}_{S \smallsetminus \{z_i\}}((\xi_i)^{\circ})$ and of the form $(x,z'_i)$ with $x \in \mathrm{Odd}_{S \smallsetminus \{z_i\}}((\xi_i)_{\dagger})$ is closed under moves and contains $\xi_i$.
Now by the above argument, for each $\zeta \in \mathcal{E}(r)$ with $\zeta_{\dagger} \in \mathrm{Odd}(r)$, the element $\zeta$ appears in one of the above subsets of $\mathcal{E}(r)$.
This implies that, for each pair $x,y$ of elements of $\mathrm{Odd}(r)$ for which $m(x,y)$ is an even integer, both $x$ and $y$ belong to the finite set $\{(\xi_i)_{\dagger},z'_i \mid 1 \leq i \leq k\}$ (apply the above argument to $(x,y) \in \mathcal{E}(r)$ and $(y,x) \in \mathcal{E}(r)$).
Hence the number of such pairs $(x,y)$ is finite, concluding the proof of Proposition \ref{prop:when_inner_generators_are_finite__acyclic}.
\end{proof}
From now, we consider the remaining case that $(\Gamma^{\mathrm{odd}})_{\sim r}$ contains a nontrivial cycle.
First, we have the following:
\begin{prop}
\label{prop:when_inner_generators_are_finite__with_cycle__if}
Let $r \in S$.
If $(\Gamma^{\mathrm{odd}})_{\sim r}$ contains a nontrivial cycle and Condition \ref{item:thm_condition_for_finite_generation_with_cycle__core} in Theorem \ref{thm:condition_for_finite_generation_with_cycle} is satisfied, then the number of inner generators is finite.
\end{prop}
\begin{proof}
The claim is obvious when Condition \ref{item:thm_condition_for_finite_generation_with_cycle__core_not_even} is satisfied.
For the remaining case that Condition \ref{item:thm_condition_for_finite_generation_with_cycle__core_A_tilde} or Condition \ref{item:thm_condition_for_finite_generation_with_cycle__core_bipyramid} is satisfied, by Remark \ref{rem:isomorphism_by_conjugation}, we may assume without loss of generality that $r \in K$.
Moreover, now we have $\mathcal{E}^{\mathrm{Odd}(r)}(r) = \mathcal{E}^K(r)$, while we have $\pi_1(\Gamma^{\mathrm{odd}};x,y) = \pi_1(\Gamma_K^{\mathrm{odd}};x,y)$ for any $x,y \in K$ by Lemma \ref{lem:fundamental_group_for_subgraph}.
Therefore, we may assume without loss of generality that $\mathrm{Odd}(r) = K$.

First, suppose that Condition \ref{item:thm_condition_for_finite_generation_with_cycle__core_A_tilde} is satisfied.
Then by Example \ref{ex:cycle_odd_length} and Example \ref{ex:cycle_even_length}, for each $\xi \in \mathcal{E}^K(r)$, there exist only a finite number of inner generators of the form $s_r(c;\xi)$.
Now the claim follows from the finiteness of $\mathcal{E}^K(r)$.

On the other hand, suppose that Condition \ref{item:thm_condition_for_finite_generation_with_cycle__core_bipyramid} is satisfied.
Then we have $\mathcal{E}^K(r) = \{(x_1,x_2),(x_2,x_1)\}$.
Now for each $i \in \{1,2\}$, by the shape of $\Gamma_K$, the group $\pi_1(\Gamma^{\mathrm{odd}};x_i)$ is generated by the elements of the form $[(x_i,y',x_{3-i},y,x_i)]$ where $y$ and $y'$ are distinct elements of $K \smallsetminus \{x_1,x_2\}$, and each generator $[(x_i,y',x_{3-i},y,x_i)]$ is realized by a sequence of E-moves $(x_i,x_{3-i}) \mapsto (x_{3-i},x_i) \mapsto (x_i,x_{3-i})$ where the first E-move is in $\{x_i,y,x_{3-i}\}$ and the second E-move is in $\{x_i,y',x_{3-i}\}$.
By the property and Theorem \ref{thm:presentation_of_orthogonal_reflection_subgroup}, all the inner generators of the form $s_r(c;(x_i,x_{3-i}))$ coincide with each other.
Hence the claim holds, concluding the proof of Proposition \ref{prop:when_inner_generators_are_finite__with_cycle__if}.
\end{proof}
Secondly, we have the following:
\begin{prop}
\label{prop:when_inner_generators_are_finite__with_cycle__only_if}
Let $r \in S$.
If $(\Gamma^{\mathrm{odd}})_{\sim r}$ contains a nontrivial cycle and the number of inner generators is finite, then Condition \ref{item:thm_condition_for_finite_generation_with_cycle__core} in Theorem \ref{thm:condition_for_finite_generation_with_cycle} is satisfied.
\end{prop}
\begin{proof}
It suffices to show that Condition \ref{item:thm_condition_for_finite_generation_with_cycle__core_A_tilde} or Condition \ref{item:thm_condition_for_finite_generation_with_cycle__core_bipyramid} in Theorem \ref{thm:condition_for_finite_generation_with_cycle} is satisfied if $m(x_1,x_2)$ is an even integer for some $x_1,x_2 \in \mathrm{Odd}(r)$.
Take such a pair $x_1,x_2$.
Then by Proposition \ref{prop:cycle_contains_even_edges} applied to $I = \mathrm{Odd}(r)$, both $x_1$ and $x_2$ are contained in every nontrivial cycle in $(\Gamma^{\mathrm{odd}})_{\sim r} = (\Gamma_{\mathrm{Odd}(r)}^{\mathrm{odd}})_{\sim r}$.
On the other hand, by Lemma \ref{lem:shape_of_cycle} applied to $I = \mathrm{Odd}(r)$, every nontrivial full cycle in $(\Gamma^{\mathrm{odd}})_{\sim r}$ is of one of the two types in the statement of Lemma \ref{lem:shape_of_cycle}.

We consider the case that a nontrivial full cycle $P$ in $(\Gamma^{\mathrm{odd}})_{\sim r}$ is of type $\widetilde{A}_{n-1}$ with $n \geq 4$.
In this case, for each vertex $x$ in $P$, there exists a vertex $y$ in $P$ for which $(x,y) \in \mathcal{E}(r)$.
Now Proposition \ref{prop:cycle_contains_even_edges} implies that $x$ is contained in every nontrivial full cycle in $(\Gamma^{\mathrm{odd}})_{\sim r}$, therefore $P$ is contained in every nontrivial full cycle in $(\Gamma^{\mathrm{odd}})_{\sim r}$.
This implies that $P$ is the unique (up to inverse and cyclic shift) nontrivial full cycle in $(\Gamma^{\mathrm{odd}})_{\sim r}$, therefore $K := V(P)$ is the core of $(\Gamma^{\mathrm{odd}})_{\sim r}$.
Hence by Proposition \ref{prop:cycle_contains_even_edges}, Condition \ref{item:thm_condition_for_finite_generation_with_cycle__core_A_tilde} is satisfied.

From now, we consider the remaining case that every nontrivial full cycle in $(\Gamma^{\mathrm{odd}})_{\sim r}$, which contains both $x_1$ and $x_2$ as shown above, is as in Figure \ref{fig:diamond}.
Note that such a nontrivial full cycle exists by the hypothesis.
Let $K$ denote the union of the vertex sets of all nontrivial full cycles in $(\Gamma^{\mathrm{odd}})_{\sim r}$.
Then $K$ is the core of $(\Gamma^{\mathrm{odd}})_{\sim r}$ and, by the property of the nontrivial full cycles mentioned above, the shape of $\Gamma_K$ is as in Condition \ref{item:thm_condition_for_finite_generation_with_cycle__core_bipyramid} in Theorem \ref{thm:condition_for_finite_generation_with_cycle}.
Now Proposition \ref{prop:cycle_contains_even_edges} implies that Condition \ref{item:thm_condition_for_finite_generation_with_cycle__core_bipyramid} in Theorem \ref{thm:condition_for_finite_generation_with_cycle} is satisfied.
Hence the proof of Proposition \ref{prop:cycle_contains_even_edges} is concluded.
\end{proof}
Now Theorem \ref{thm:condition_for_finite_generation_with_cycle} follows by combining Proposition \ref{prop:when_outer_generators_are_finite}, Proposition \ref{prop:when_inner_generators_are_finite__with_cycle__if} and Proposition \ref{prop:when_inner_generators_are_finite__with_cycle__only_if}.
On the other hand, Theorem \ref{thm:condition_for_finite_generation_without_cycle} follows by combining Proposition \ref{prop:when_outer_generators_are_finite} and Proposition \ref{prop:when_inner_generators_are_finite__acyclic}.
Hence the proofs of Theorem \ref{thm:condition_for_finite_generation_with_cycle} and Theorem \ref{thm:condition_for_finite_generation_without_cycle} are concluded.

\end{document}